\documentclass[final]{amsart}
\usepackage[dvips,final]{graphicx}
\usepackage{enumerate}

\usepackage[ps2pdf,draft=false]{hyperref}

\DeclareMathOperator{\Li}{Li}
\DeclareMathOperator{\CAT}{CAT}

\theoremstyle{plain}
\newtheorem{theorem}{Theorem}[section]
\newtheorem{lemma}[theorem]{Lemma}
\newtheorem{proposition}[theorem]{Proposition}

\theoremstyle{definition}
\newtheorem{definition}[theorem]{Definition}

\theoremstyle{remark}

\newtheorem{remark}[theorem]{Remark}

\numberwithin{equation}{section}
\numberwithin{figure}{section}
\numberwithin{table}{section}

\newcommand{\E}[1]{\mathbf{E}\,#1}
\newcommand{\Var}[1]{\mathbf{Var}\,#1}

\begin{document}
\title{Limiting distributions for additive functionals on Catalan
trees}
\thanks{Research for both authors supported by NSF grants DMS-9803780 and DMS-0104167,
and by The Johns Hopkins University's Acheson J.~Duncan Fund
for the Advancement of Research in Statistics.  Research for the second author supported
by NSF grant 0049092 and carried out primarily while this author was affiliated with what is now
the Department of Applied Mathematics and Statistics at The Johns Hopkins University.}
\author{James Allen Fill}
\email[James Allen Fill]{jimfill@jhu.edu}
\urladdr[James Allen Fill]{http://www.ams.jhu.edu/\~{}fill/}
\author{Nevin Kapur}
\email[Nevin Kapur]{nkapur@cs.caltech.edu}
\urladdr[Nevin Kapur]{http://www.cs.caltech.edu/\~{}nkapur/}
\date{June 4, 2003.  Revised April 1, 2004.}
\address[James Allen Fill]{Department of Applied Mathematics and
  Statistics, The  Johns Hopkins 
  University, 3400 N. Charles St., Baltimore MD 21218, USA}
\address[Nevin Kapur]{Department of Computer Science, California
  Institute of Technology, MC 256-80, 1200 E.~California Blvd.,
  Pasadena CA 91125, USA}

\begin{abstract}
  Additive tree functionals represent the cost of many
  divide-and-conquer algorithms.  We derive the limiting distribution
  of the additive functionals induced by toll functions of the form
  (a)~\(n^\alpha\) when \(\alpha > 0\) and (b)~\(\log{n}\) (the
  so-called shape functional) on uniformly distributed binary
  trees, sometimes called Catalan trees.  The Gaussian law obtained in
  the latter case complements the central limit theorem for the shape
  functional under the random permutation model.
  Our results give rise to an
  apparently new family of distributions containing the Airy
  distribution ($\alpha = 1$) and the normal distribution [case~(b),
  and case~(a) as $\alpha \downarrow 0$].  The main theoretical tools
  employed are recent results relating asymptotics of the generating
  functions of sequences to those of their Hadamard product, and the
  method of moments.
\end{abstract}

\maketitle

\section{Introduction}
\label{sec:introduction}

Binary
trees are fundamental data structures in computer
science, with primary application in searching and sorting.  For
background we refer the reader to Chapter~2 of the excellent
book~\cite{MR93f:68045}.  In this article we consider additive
functionals defined on uniformly distributed binary
trees (sometimes called Catalan trees) induced by two types of toll
sequences [$(n^{\alpha})$ and $(\log n)$].  (See the simple
Definition~\ref{def:additive-functional}.)  Our main results,
Theorems~\ref{thm:limit-dist} and~\ref{thm:shape_clt}, establish the
limiting distribution for these induced functionals.

A competing model of randomness for binary
trees---one used for binary search trees---is the
\emph{random permutation model} (RPM); see Section~2.3
of~\cite{MR93f:68045}.  While there has been much study of
additive functionals under the RPM (see, for
example,~\cite[Section~3.3]{MR93f:68045}
and~\cite{MR92f:68028,MR97f:68021,01808498,MR1871558}), little
attention has been paid to the distribution of functionals defined on
binary
trees under the uniform (Catalan) model of randomness.
Fill~\cite{MR97f:68021} argued that the functional corresponding to
the toll sequence \( (\log{n}) \) serves as a crude measure of the
``shape'' of a binary
tree, and explained how this functional
arises in connection with the move-to-root self-organizing scheme for
dynamic maintenance of binary search trees.  He derived a central
limit theorem under the RPM, but obtained only asymptotic information
about the mean and variance under the Catalan model.  (The latter
results were rederived in the extension~\cite{MR99j:05171} from binary
trees to simply generated rooted trees.)  In this paper
(Theorem~\ref{thm:shape_clt}) we show that there is again asymptotic
normality under the Catalan model.

In~\cite[Prop.~2]{MR88h:68033} Flajolet and Steyaert gave
order-of-growth information about the mean of functionals induced by
tolls of the form \( n^\alpha \).  (The motivation is to build a
``repertoire'' of tolls from which the behavior of more complicated
tolls can be deduced by combining elements from the repertoire.  The
corresponding results under the random permutation model were derived
by Neininger~\cite{MR1910527}.)
Tak{\'a}cs established the limiting (Airy) distribution of path length
in Catalan trees~\cite{MR92m:60057,MR93e:60175,MR92k:60164}, which is
the additive functional for the toll~$n-1$.  The additive functional
for the toll $n^2$ arises in the study of the Wiener index of the tree
and has been analyzed by Janson~\cite{janson:_wiener}. In this paper
(Theorem~\ref{thm:limit-dist}) we obtain the limiting distribution for
Catalan trees for toll $n^{\alpha}$ for any $\alpha > 0$.  The family
of limiting distributions appears to be new.  In most cases we have a
description of the distribution only in terms of its moments, although
other descriptions in terms of Brownian excursion, as for the Airy
distribution and the limiting distribution for the Wiener index, may
be possible.  This is currently under investigation by the authors in
collaboration with others.

The uniform model on binary trees has also been used
recently by Janson~\cite{janson02:_ideal} in the analysis of an
algorithm of Koda and Ruskey~\cite{MR94f:94015} for listing ideals in
a forest poset.

This paper serves as the first example of the application of recent
results~\cite{FFK}, extending singularity analysis~\cite{MR90m:05012},
to obtain limiting distributions.  In~\cite{FFK}, it is shown how the
asymptotics of generating functions of sequences relate to those of
their Hadamard product.  First moments for our problems were treated
in~\cite{FFK} and a sketch of the technique we employ was presented
there. (Our approach to obtaining asymptotics of Hadamard products of
generating functions differs only marginally from the Zigzag Algorithm
as presented in~\cite{FFK}.)  As will be evident soon, Hadamard
products occur naturally when one is analyzing moments of additive
tree functionals.  The program we carry out allows a fairly mechanical
derivation of the asymptotics of moments of each order, thereby
facilitating application of the method of moments.  Indeed,
preliminary investigations suggest that the techniques we develop are
likewise applicable to the wider class of simply generated trees; this
is work in progress.

The organization of this paper is as follows.
Section~\ref{sec:preliminaries} establishes notation and states
certain preliminaries that will be used in the subsequent proofs.  In
Section~\ref{sec:toll-sequence-nalpha} we consider the toll sequence
$(n^\alpha)$ for general $\alpha > 0$.  In
Section~\ref{sec:asympotics-mean} we compute the asymptotics of the
mean of the corresponding additive functional.  In
Section~\ref{sec:higher-moments} the analysis diverges slightly as the
nature of asymptotics of the higher moments differs depending on the
value of~\(\alpha\).  Section~\ref{sec:asymptotics-moments} employs
singularity analysis~\cite{MR90m:05012} to derive the asymptotics of
moments of each order.  In Section~\ref{sec:limit-distr} we use the
results of Section~\ref{sec:asymptotics-moments} and the method of
moments to derive the limiting distribution of the additive tree
functional.  In Section~\ref{sec:shape-functional} we employ the
approach again to obtain a normal limit theorem for the shape
functional.  Finally, in Section~\ref{sec:suff-cond-asympt}, we
present heuristic arguments that may lead to the identification of
toll sequences giving rise to a normal limit.


\newpage

\section{Notation and Preliminaries}
\label{sec:preliminaries}

\subsection{Additive tree functionals}
\label{sec:addit-tree-funct}

We first establish some notation.  Let~\(T\) be a binary
tree.
We use~\(|T|\) to denote the number of
nodes
in~\(T\). Let $L(T)$ and
$R(T)$ denote, respectively, the left and right subtrees rooted at the
children of the root of $T$.
\begin{definition}
  \label{def:additive-functional}
  A functional~\( f \) on binary
  trees is called an \emph{additive tree functional} if it
  satisfies the recurrence
  \begin{equation*}
    \label{eq:2.1}
    f(T) = f( L(T) ) + f( R(T) )+ b_{|T|},
  \end{equation*}
  for any tree~\( T \) with~\( |T| \geq 1 \). Here \( (b_n)_{n \geq
  1} \) is a given sequence, henceforth called the \emph{toll
  function}. 
\end{definition}

We analyze 
additive functionals
defined on
binary
trees uniformly distributed over $\{ T\!:\, |T|=n \}$ for
given $n$.  Let $X_n$ be such an additive
functional induced by the toll sequence $(b_n)$.
It is well known that the number of binary
trees on $n$ nodes
is counted by the $n$th Catalan number
\[
\beta_n := \frac1{n+1}\binom{2n}{n},
\]
with generating function
\begin{equation*}
  \label{eq:36}
  \CAT(z) := \sum_{n=0}^\infty \beta_n z^n = \frac1{2z}(1 - \sqrt{1-4z}).
\end{equation*}
In our subsequent analysis we will make use of the identity
\begin{equation}
  \label{eq:50}
  z \CAT^2(z) = \CAT(z) - 1.
\end{equation}
The mean of the cost function \(a_n := \E{X_n} \) can be
obtained recursively by conditioning on
the size of $L(T)$
as
\begin{equation*}
  a_n = \sum_{j=1}^n \frac{\beta_{j-1}\beta_{n-j}}{\beta_n}(a_{j-1} +
  a_{n-j}) + b_n, \qquad n \geq 1.
\end{equation*}
This recurrence can be rewritten as
\begin{equation}
  \label{eq:35}
  (\beta_n a_n) = 2 \sum_{j=1}^n (\beta_{j-1} a_{j-1}) \beta_{n-j} +
  (\beta_n b_n), \qquad n \geq 1.
\end{equation}
Recall that the \emph{Hadamard product} of two power series $F$ and
$G$, denoted by $F(z) \odot G(z)$, is the power series defined by
\begin{equation*}
  ( F \odot G)(z) \equiv F(z) \odot G(z) := \sum_{n} f_n g_n z^n,
\end{equation*}
where
\begin{equation*}
  F(z) = \sum_{n} f_n z^n \qquad \text{ and } \qquad G(z) = \sum_{n}
  g_n z^n.
\end{equation*}
Multiplying~(\ref{eq:35}) by $z^n/4^n$ and summing over $n \geq 1$ we get
\begin{equation}
  \label{eq:1}
A(z)\odot\CAT(z/4) = \frac{B(z)\odot\CAT(z/4)}{\sqrt{1-z}}, 
\end{equation}
where \(A(z)\) and \(B(z)\) are the ordinary generating functions of
\((a_n)\) and \((b_n)\) respectively.
\begin{remark}
  \label{rem:Catalan}
  Catalan numbers are ubiquitous in combinatorial applications;
  see~\cite{MR2000k:05026} for a list of 66~instances and
  \url{http://www-math.mit.edu/~rstan/ec/} for more.
\end{remark}

In the sequel the notation [\(\cdots\)] is used both for Iverson's
convention~\cite[1.2.3(16)]{knuth97}
and for the coefficient of
certain terms in the succeeding expression.  The interpretation will
be clear from the context.  For example, \( [ \alpha > 0 ] \) has the
value 1 when \( \alpha > 0 \) and the value 0 otherwise.  In contrast, \( [z^n]
F(z) \) denotes the coefficient of \( z^n \) in the series expansion
of \( F(z) \).  Throughout this paper $\Gamma$ and $\zeta$ denote
Euler's gamma function and the Riemann zeta function, respectively.

\subsection{Singularity analysis}
\label{sec:singularity-analysis}

\emph{Singularity analysis} is a systematic complex-analytic technique
that relates asymptotics of sequences to singularities of their
generating functions.  The applicability of singularity analysis rests
on the technical condition of~\emph{$\Delta$-regularity}.  Here is the
definition.  See~\cite{FFK} or~\cite{MR90m:05012} for further background.
\begin{definition}\label{def:delta-regular}
  A function defined by a Taylor series about the origin with radius
  of convergence equal to~$1$ is \emph{$\Delta$-regular} if
  it can be analytically continued in a domain
\begin{equation*}
   \Delta(\phi,\eta) := \{z: |z| < 1 + \eta, |\arg(z-1)| > \phi\},
\end{equation*}
for some $\eta > 0$ and $0 < \phi < \pi/2$. A function $f$ is said to
admit a \emph{singular expansion} at $z=1$ if it is $\Delta$-regular
and
\begin{equation*}
   f(z) = \sum_{j=0}^J c_j(1-z)^{\alpha_j} + O(|1-z|^A)
\end{equation*}
uniformly in $z \in \Delta(\phi,\eta)$,
for a sequence of complex numbers
$(c_j)_{0 \leq j \leq J}$ and an increasing sequence of real numbers
$(\alpha_j)_{0 \leq j \leq J}$ satisfying $\alpha_j < A$.  It is said
to satisfy a singular expansion \emph{``with logarithmic terms''} if,
similarly, 
\begin{equation*}
f(z) = \sum_{j=0}^J c_j\left(L(z)\right)(1-z)^{\alpha_j} + O(|1-z|^A),
\qquad
L(z):=\log\frac{1}{1-z},
\end{equation*}
where~each $c_j(\cdot)$ is a polynomial.
\end{definition}
Following established terminology, when a function has a singular
expansion with logarithmic terms we shall say that it is \emph{amenable} to singularity
analysis.

Recall the definition of the \emph{generalized polylogarithm}:
\begin{definition}\label{def:li}
  For $\alpha$ an arbitrary complex number and $r$ a nonnegative
  integer, the \emph{generalized polylogarithm} function
  $\Li_{\alpha,r}$ is defined  for $|z| < 1$ by
  \begin{equation*}
    \label{eq:4.1.50}
    \Li_{\alpha,r}(z) := \sum_{n=1}^\infty \frac{(\log{n})^r}{n^\alpha} z^n.
  \end{equation*}
\end{definition}
The key property of the generalized polylogarithm that we will employ
is
\begin{equation*}
  \Li_{\alpha,r} \odot \Li_{\beta,s} = \Li_{\alpha+\beta,r+s}.
\end{equation*}
We will also make extensive use of the following consequences of
the singular expansion of the generalized polylogarithm.
Neither this lemma nor the ones following make any claims about
uniformity in $\alpha$ or $r$.
Note that $\Li_{1,0}(z) = L(z) = \log\bigl((1-z)^{-1}\bigr)$.
\begin{lemma}
  \label{lem:litoomz}
  For any real $\alpha < 1$ and nonnegative integer~$r$, we have the
  singular expansion
  \[
  \Li_{\alpha,r}(z) = \sum_{k=0}^r \lambda_k^{(\alpha,r)}
  (1-z)^{\alpha-1} L^{r-k}(z) +
  O(|1-z|^{\alpha-\epsilon}) + (-1)^r \zeta^{(r)}(\alpha)[\alpha > 0],
  \]
  where $\lambda_k^{(\alpha,r)} \equiv \binom{r}{k}
  \Gamma^{(k)}(1-\alpha)$ and $\epsilon > 0$ is arbitrarily small.
\end{lemma}
\begin{proof}
  By Theorem~1
  in~\cite{MR2000a:05015},
  \begin{equation}\label{eq:37}
    \Li_{\alpha,0}(z) \sim \Gamma(1-\alpha) t^{\alpha-1} + \sum_{j
    \geq 0} \frac{(-1)^j}{j!} \zeta(\alpha-j) t^j, \qquad t = -\log z
    = \sum_{l=1}^\infty \frac{(1-z)^l}{l},
  \end{equation}
  and for any positive integer~$r$,
  \begin{equation*}
    \Li_{\alpha,r}(z) = (-1)^r \frac{\partial^r}{\partial\alpha^r}
    \Li_{\alpha,0}(z).
  \end{equation*}
  Moreover, as also shown in~\cite{MR2000a:05015}, the singular
  expansion for $\Li_{\alpha,r}$ is obtained by
  performing the indicated differentiation of~(\ref{eq:37}) term-by-term.
  To establish the claim we set \( f = \Gamma(1-\alpha) \) and \( g =
  t^{\alpha-1} \) in the general formula for the $r$th derivative of a
  product:
   \begin{equation*}
    (fg)^{(r)} = \sum_{k=0}^r \binom{r}{k} f^{(k)} g^{(r-k)}
   \end{equation*}
   to first obtain
   \begin{equation*}
     (-1)^r \frac{\partial^r}{\partial\alpha^r}
    [\Gamma(1-\alpha)t^{\alpha-1}] = (-1)^r \sum_{k=0}^r \binom{r}{k}
    (-1)^k\Gamma^{(k)}(1-\alpha) t^{\alpha-1} (\log t)^{r-k}
   \end{equation*}
   The claim then follows easily.
\end{proof}
The following ``inverse'' of Lemma~\ref{lem:litoomz} is very useful
for computing with Hadamard products.
\begin{lemma}
  \label{lem:omztoli}
  For any real $\alpha < 1$ and nonnegative integer~$r$, there exists
  a region~$\Delta(\phi,\eta)$ as in Defintion~\ref{def:delta-regular}
  such that
  \[
  (1-z)^{\alpha-1}L^r(z) = \sum_{k=0}^r
  \mu_k^{(\alpha,r)} \Li_{\alpha,r-k}(z) + O(|1-z|^{\alpha-\epsilon})
  + c_r(\alpha)[\alpha > 0]
  \]
  holds uniformly in $z \in \Delta(\phi,\eta)$,
  where $\mu_0^{(\alpha,r)}=1/\Gamma(1-\alpha)$, \( c_r(\alpha) \) is
  a constant, and
  $\epsilon > 0$ is  arbitrarily small.
\end{lemma}
\begin{proof}
  We use induction on~$r$. For $r=0$ we have
  \[
  \Li_{\alpha,0}(z) = \Gamma(1-\alpha)(1-z)^{\alpha-1} +
  O(|1-z|^{\alpha-\epsilon}) + \zeta(\alpha) [ \alpha > 0 ]
  \]
  and the claim is verified with
  \begin{equation*}
    \mu_0^{(\alpha,0)} = \frac{1}{\Gamma(1-\alpha)} \qquad  \text{ and
    } \qquad c_0(\alpha) = -\frac{\zeta(\alpha)}{\Gamma(1-\alpha)}.
  \end{equation*}
  Let $r \geq 1$.  Then using
  Lemma~\ref{lem:litoomz} and the induction hypothesis we get
  \begin{align*}
    &\Li_{\alpha,r}(z) \\&=
    \Gamma(1-\alpha)(1-z)^{\alpha-1}
    L^r(z) \\
    & \quad + \sum_{k=1}^r \lambda_{k}^{(\alpha,r)} \left[
      \sum_{l=0}^{r-k}\mu_l^{(\alpha,r-k)} \Li_{\alpha,r-k-l}(z) +
    O(|1-z|^{\alpha-\epsilon}) +
      c_{r-k}(\alpha) [\alpha > 0 ] \right]\\
    & \quad + O(|1-z|^{\alpha-\epsilon})
    + (-1)^r\zeta^{(r)}(\alpha)[\alpha > 0]\\
    &=
    \Gamma(1-\alpha)(1-z)^{\alpha-1}L^r(z)
    + \sum_{k=1}^r \lambda_{k}^{(\alpha,r)} \sum_{s=0}^{r-k}
    \mu_{r-k-s}^{(\alpha,r-k)} \Li_{\alpha,s}(z)  \\
    & \quad {}+  O(|1-z|^{\alpha-\epsilon}) +
    \left(
      \sum_{k=1}^r \lambda_k^{(\alpha,r)} c_{r-k}(\alpha)
      + (-1)^r\zeta^{(r)}(\alpha)
    \right)[\alpha > 0]\\
    &=
    \Gamma(1-\alpha)(1-z)^{\alpha-1}
    L^r(z)
    + \sum_{s=0}^{r-1} \nu_s^{(\alpha,r)} \Li_{\alpha,s}(z)\\
    & \quad +  O(|1-z|^{\alpha-\epsilon}) + \gamma_r(\alpha)[\alpha >
    0],
  \end{align*}
  where, for $0 \leq s \leq r-1$,
  \[
  \nu_s^{(\alpha,r)} := \sum_{k=1}^{r-s} \lambda_k^{(\alpha,r)}
  \mu_{r-s-k}^{(\alpha,r-k)},
  \]
  and where
  \begin{equation*}
    \gamma_r(\alpha) := \sum_{k=1}^r
    \lambda_k^{(\alpha,r)}c_{r-k}(\alpha) + (-1)^r \zeta^{(r)}(\alpha).
  \end{equation*}
  Setting
  \begin{equation*}
    \mu_0^{(\alpha,r)} = \frac{1}{\Gamma(1-\alpha)},\qquad 
    \mu_k^{(\alpha,r)} =
    -\frac{\nu_{r-k}^{(\alpha,r)}}{\Gamma(1-\alpha)}, \quad 1 \leq k
    \leq r,
  \end{equation*}
  and
  \begin{equation*}
    c_r(\alpha) = - \frac{\gamma_r(\alpha)}{\Gamma(1-\alpha)},
  \end{equation*}
  the result follows.
\end{proof}

For the calculation of the mean, the following refinement of a special
case of Lemma~\ref{lem:litoomz} is required.  It is a simple
consequence of Theorem~1 of~\cite{MR2000a:05015}.
\begin{lemma}
  \label{lem:li0omz}
  When \(\alpha < 0\), we have the singular expansion
  \[
  \Li_{\alpha,0}(z) = \Gamma(1-\alpha)(1-z)^{\alpha-1} -
  \Gamma(1-\alpha)\frac{1-\alpha}2 (1-z)^{\alpha} +
  O(|1-z|^{\alpha+1}) + \zeta(\alpha)[\alpha > -1] .
  \]
\end{lemma}
For the sake of completeness, we state a result of particular relevance
from~\cite{FFK}.
\begin{theorem}
  \label{thm:hadamard}
  If $f$ and $g$ are amenable to singularity analysis and
  \begin{equation*}
    f(z) = O(|1-z|^a) \qquad\text{ and }\qquad g(z) = O(|1-z|^b)
  \end{equation*}
  as $z \to 1$,  then  $f \odot g$ is also amenable to singularity
  analysis.  Furthermore
  \begin{enumerate}[(a)]
  \item If $a+b+1 < 0$ then
    \begin{equation*}
      f(z) \odot g(z) = O(|1-z|^{a+b+1}).
    \end{equation*}\label{item:hadamard1}
  \item If $k < a+b+1 < k+1$ for some integer $-1 \leq k < \infty$,
    then
    \begin{equation*}
      f(z) \odot g(z) = \sum_{j=0}^k \frac{(-1)^j}{j!} (f \odot
      g)^{(j)}(1) (1-z)^j + O(|1-z|^{a+b+1}).
    \end{equation*}
    \label{item:hadamard2}
  \item   If $a+b+1$ is a nonnegative integer then
    \begin{equation*}
      f(z) \odot g(z) = \sum_{j=0}^{a+b} \frac{(-1)^j}{j!} (f \odot
      g)^{(j)}(1) (1-z)^j + O(|1-z|^{a+b+1}|L(z)|).
    \end{equation*}    \label{item:hadamard3}
  \end{enumerate}
\end{theorem}

\section{\texorpdfstring{The toll sequence $(n^\alpha)$}{The toll
    sequence n^\alpha}}
\label{sec:toll-sequence-nalpha}

In this section we consider additive functionals when the toll
function $b_n$ is $n^\alpha$ with $\alpha > 0$.

\subsection{Asympotics of the mean}
\label{sec:asympotics-mean}
The main result of this Section~\ref{sec:asympotics-mean} is a
singular expansion for $A(z) \odot \CAT(z/4)$.  The result
is~\eqref{eq:9},~\eqref{eq:7}, or~\eqref{eq:8} according as $\alpha <
1/2$, $\alpha = 1/2$, or $\alpha > 1/2$.

Since $b_n = n^\alpha$, by definition \(B=\Li_{-\alpha,0}\). Thus,
by Lemma~\ref{lem:li0omz},
\[
B(z) = \Gamma(1+\alpha)(1-z)^{-\alpha-1} -
\Gamma(1+\alpha)\frac{\alpha+1}2(1-z)^{-\alpha}  +
O(|1-z|^{-\alpha+1}) + \zeta(-\alpha)[\alpha < 1].
\]
We will now use~(\ref{eq:1}) to obtain the asymptotics of the mean.

First we treat the case
\(\alpha < 1/2\).
From the singular expansion \(\CAT(z/4) = 2 + O(|1-z|^{1/2})\) as \(z
\to 1\), we have, by part~(\ref{item:hadamard2}) of
Theorem~\ref{thm:hadamard},
\[
B(z)\odot\CAT(z/4)  = C_0 + O(|1-z|^{-\alpha+\tfrac12}),
\]
where
\begin{equation*}
  \label{eq:20}
  C_0 := B(z)\odot\CAT(z/4) \Bigr\rvert_{z=1} = \sum_{n=1}^\infty n^\alpha
  \frac{\beta_n}{4^n}.
\end{equation*}
We now already know the constant term in the singular
expansion of \(B(z)\odot\CAT(z/4)\) at \(z=1\) and henceforth we  need
only compute
lower-order terms.  The constant \(\bar{c}\)
is used in the sequel to denote an unspecified (possibly~0) constant,
possibly different at each appearance.

Let's write \(B(z) = L_1(z) + R_1(z)\), and \(\CAT(z/4) = L_2(z) +
R_2(z)\), where
\begin{align*}
  L_1(z) &:= \Gamma(1+\alpha)(1-z)^{-\alpha-1} -
  \Gamma(1+\alpha)\frac{\alpha+1}2(1-z)^{-\alpha} + \zeta(-\alpha),\\
  R_1(z) &:= B(z) - L_1(z) = O(|1-z|^{1-\alpha}),\\
  L_2(z) &:= 2(1 - (1-z)^{1/2}),\\
  R_2(z) &:= \CAT(z/4) - L_2(z) = O(|1-z|).
\end{align*}
We will analyze each of the four Hadamard products separately. First,
\begin{align*}
  L_1(z)\odot{}L_2(z) &= -
 2\Gamma(1+\alpha)\bigl[(1-z)^{-\alpha-1}\odot(1-z)^{1/2}\bigr]\\
 & \quad +  2\Gamma(1+\alpha)\frac{\alpha+1}2
 \bigl[(1-z)^{-\alpha}\odot(1-z)^{1/2}\bigr] + \bar{c}.
\end{align*}
By Theorem~4.1 of~\cite{FFK},
\[
(1-z)^{-\alpha-1}\odot(1-z)^{1/2} = \bar{c} +
\frac{\Gamma(\alpha-\tfrac12)}{\Gamma(\alpha+1)\Gamma(-1/2)}
(1-z)^{-\alpha+\tfrac12} + O(|1-z|),
\]
and
\[
(1-z)^{-\alpha}\odot(1-z)^{1/2} = \bar{c} + O(|1-z|)
\]
by another application of part~(\ref{item:hadamard2}) of
Theorem~\ref{thm:hadamard}, this time with $k=1$. Hence
\[
L_1(z)\odot{}L_2(z) = \bigl[L_1(z)\odot{}L_2(z)\bigr]\Bigr\rvert_{z=1} +
\frac{\Gamma(\alpha-\tfrac12)}{\sqrt\pi}
(1-z)^{-\alpha+\tfrac12} + O(|1-z|).
\]
The other three Hadamard products are easily handled as
\begin{align*}
  L_1(z)\odot{}R_2(z) &= \bigl[L_1(z)\odot{}R_2(z)\bigr]\Bigr\rvert_{z=1} +
  O(|1-z|^{-\alpha+1}), \\
  L_2(z)\odot{}R_1(z) &= \bigl[L_2(z)\odot{}R_1(z)\bigr]\Bigr\rvert_{z=1} +
  O(|1-z|), \\
  R_1(z)\odot{}R_2(z) &= \bigl[R_1(z)\odot{}R_2(z)\bigr]\Bigr\rvert_{z=1} +
  O(|1-z|).
\end{align*}
Putting everything together, we get
\[ B(z)\odot\CAT(z/4) = C_0 +
\frac{\Gamma(\alpha-\tfrac12)}{\sqrt\pi}
(1-z)^{-\alpha+\tfrac12} + O(|1-z|^{-\alpha+1}).
\]
Using this in~\eqref{eq:1}, we get
\begin{equation}
  \label{eq:9}
  A(z)\odot\CAT(z/4) = C_0(1-z)^{-1/2} +
  \frac{\Gamma(\alpha-\tfrac12)}{\sqrt\pi} (1-z)^{-\alpha} +
  O(|1-z|^{-\alpha+\tfrac12}).
\end{equation}

To treat the case \(\alpha \geq 1/2\) we make use of the estimate
\begin{equation}
  \label{eq:21}
  (1-z)^{1/2} = \frac1{\Gamma(-1/2)}[ \Li_{3/2,0}(z) - \zeta(3/2) ]
  + O(|1-z|),
\end{equation}
a consequence of Theorem~1 of~\cite{MR2000a:05015},
so that
\[
B(z) \odot (1-z)^{1/2} = \Li_{-\alpha,0}(z) \odot (1-z)^{1/2} =
\frac{1}{\Gamma(-1/2)} \Li_{\tfrac32-\alpha,0}(z) + R(z),
\]
where
\begin{equation}
  \label{eq:55}
  R(z) =
  \begin{cases}
    \bar{c} + O(|1-z|^{1-\alpha}) & 1/2 \leq \alpha < 1 \\
    O(|L(z)|)                     & \alpha=1            \\
    O(|1-z|^{1-\alpha})           & \alpha > 1.
  \end{cases}
\end{equation}
Hence
\[
B(z)\odot\CAT(z/4) = -\frac2{\Gamma(-1/2)}\Li_{\tfrac32-\alpha,0}(z) +
\widetilde{R}(z), 
\]
where $\widetilde{R}$, like $R$, satisfies~\eqref{eq:55} (with
a possibly different $\bar{c}$).  When \(\alpha=1/2\), this gives us
\[
B(z)\odot\CAT(z/4) = -\frac2{\Gamma(-1/2)}L(z) +
\bar{c} + O(|1-z|^{1/2}),
\]
so that
\begin{equation}
  \label{eq:7}
  A(z)\odot\CAT(z/4) =
  \frac1{\sqrt\pi}(1-z)^{-1/2}L(z) +
  \bar{c}(1-z)^{-1/2} + O(1).
\end{equation}
For \(\alpha > 1/2\) another singular expansion leads to the conclusion
that
\begin{equation}
  \label{eq:8}
  A(z)\odot\CAT(z/4) =
  \frac{\Gamma(\alpha-\tfrac12)}{\sqrt\pi}(1-z)^{-\alpha} + \widehat{R}(z),
\end{equation}
where
\begin{equation*}
  \widehat{R}(z) = 
  \begin{cases}
    O(|1-z|^{-\tfrac12})        & 1/2 < \alpha < 1 \\
    O(|1-z|^{-\tfrac12}|L(z)|)  & \alpha = 1       \\
    O(|1-z|^{-\alpha+\tfrac12}) & \alpha > 1.
  \end{cases}
\end{equation*}

We defer deriving the asymptotics of \(a_n\) until
Sections~\ref{sec:higher-moments}--\ref{sec:asymptotics-moments}.

\subsection{Higher moments}
\label{sec:higher-moments}

We will analyze separately the cases \(0 < \alpha < 1/2\),
\mbox{$\alpha=1/2$}, and \mbox{$\alpha > 1/2$}.  The reason for this
will become evident soon; though the technique used to derive the
asymptotics is induction in each case, the induction hypothesis is
different for each of these cases.

\subsubsection{Small toll functions\texorpdfstring{ ($0 < \alpha <
    1/2$)}{}}
\label{sec:small-toll-functions}

We start by restricting ourselves to tolls of the form \(n^\alpha\)
where \(0 < \alpha < 1/2\).  In this case we observe that by
singularity analysis applied to~(\ref{eq:9}),
\begin{equation*}
  \frac{a_n \beta_n}{4^n} = \frac{C_0}{\sqrt\pi} n^{-1/2} + O(n^{-3/2}) +
  O(n^{\alpha-1}) = \frac{C_0}{\sqrt\pi} n^{-1/2} +
  O(n^{\alpha-1}),
\end{equation*}
so
\begin{equation*}
  a_n = n^{\tfrac32}[1 + O(n^{-1})][ C_0n^{-\tfrac12} +
  O(n^{\alpha-1})] = C_0 n + O(n^{\alpha+\tfrac12}) = (C_0 + o(1)) (n+1).
\end{equation*}
The lead-order term of the mean $a_n = \E{X_n}$
is thus linear, irrespective of the value of \(0 < \alpha < 1/2\)
(though the coefficient~$C_0$ does depend on~\(\alpha\)).  We next
perform an approximate centering to get to further dependence on \(\alpha\).


Define \(\widetilde{X}_n := X_n -
C_0(n+1)\), with \(X_0 := 0\); \(\tilde{\mu}_n(k) :=
\E{\widetilde{X}_n^k}\), with $\tilde{\mu}_n(0) = 1$ for all $n \geq
0$; and \( \hat{\mu}_n(k) :=
\beta_n\tilde{\mu}_n(k)/4^n \).  Let \(\widehat{M}_k(z)\) denote the ordinary
generating function of \(\hat{\mu}_n(k)\) in the argument $n$.

By an argument similar to the one that led to~(\ref{eq:35}), we get,
for \(k \geq 2\),
\[
\hat{\mu}_n(k) = \frac12 \sum_{j=1}^n \frac{\beta_{n-j}}{4^{n-j}}
\hat{\mu}_{j-1}(k) + \hat{r}_n(k), \qquad n \geq 1,
\]
where
\begin{align*}
\hat{r}_n(k) &:= \frac14\sum_{j=1}^n \sum_{\substack{k_1+k_2+k_3=k\\
    k_1,k_2 < k}} \binom{k}{k_1,k_2,k_3}
    \hat{\mu}_{j-1}(k_1)
    \hat{\mu}_{n-j}(k_2) b_n^{k_3}\\
    &= \frac14 \sum_{\substack{k_1+k_2+k_3=k\\
        k_1,k_2 < k}} \binom{k}{k_1,k_2,k_3} b_n^{k_3} \sum_{j=1}^n
    \hat{\mu}_{j-1}(k_1) \hat{\mu}_{n-j}(k_2),
\end{align*}
for \(n \geq 1\) and \(\hat{r}_0(k) := \hat{\mu}_0(k) = \tilde{\mu}_0(k) =
(-1)^kC_0^k\). Let \(\widehat{R}_k(z)\) denote the ordinary generating
function of \(\hat{r}_n(k)\) in the argument $n$.  Then,
mimicking~(\ref{eq:1}),
\begin{equation}
  \label{eq:4}
\widehat{M}_k(z) = \frac{\widehat{R}_k(z)}{\sqrt{1-z}}  
\end{equation}
with
\begin{equation}
  \label{eq:3}
\widehat{R}_k(z) = (-1)^kC_0^k + \sum_{\substack{k_1+k_2+k_3=k\\
    k_1,k_2<k}} \binom{k}{k_1,k_2,k_3} \bigl( B(z)^{\odot k_3}
\bigr) \odot
\bigl[\frac{z}4 \widehat{M}_{k_1}(z)\widehat{M}_{k_2}(z)\bigr],
\end{equation}
where for \(k\) a nonnegative integer
\[
B(z)^{\odot k} := \underbrace{B(z)\odot\cdots\odot B(z)}_{k}.
\]
Note that \(\widehat{M}_0(z) = \CAT(z/4)\).
\begin{proposition}
  \label{thm:0alpha14}
    Let $\epsilon > 0$ be arbitrary, and define
  \begin{equation*}
    c :=
    \begin{cases}
      2\alpha - \epsilon      & 0 < \alpha \leq  1/4 \\
      1/2 & 1/4 < \alpha < 1/2.
    \end{cases}
  \end{equation*}
  Then we have the singular expansion
  \[
  \widehat{M}_k(z) = C_k(1-z)^{-k(\alpha+\tfrac12) +\tfrac12} +
  O(|1-z|^{-k(\alpha+\tfrac{1}{2}) + \tfrac12 + c}),
  \]
  The \(C_k\)'s here are defined by the recurrence
  \begin{equation}
    \label{eq:10}
    C_k = \frac14 \sum_{j=1}^{k-1} \binom{k}{j} C_{j} C_{k-j} +
    kC_{k-1} \frac{ \Gamma(k\alpha+\tfrac{k}2-1)}{\Gamma((k-1)\alpha +
      \tfrac{k}2 -1)}, \quad k \geq 2; \quad C_1 =
    \frac{\Gamma(\alpha-\tfrac12)}{\sqrt\pi}.
  \end{equation}
\end{proposition}
\begin{proof}
  For \(k=1\) the claim is true as shown in~\eqref{eq:9} with \(C_1\) as
  defined in~\eqref{eq:10}.  We will now analyze each term
  in~\eqref{eq:3} for \(k \geq 2\).

  One can analyze separately the cases $0 < \alpha \leq 1/4$ and $1/4
  < \alpha < 1/2$.  The proof technique in either case is induction.
  We shall treat here the case $0 < \alpha \leq 1/4$; the details in the
  other case can be found in~\cite{FK-catalan-arXiv}.

  For notational convenience, define
  \(\alpha' := \alpha+\tfrac12\).  Also, observe that
  \[
  B(z)^{\odot{}k} = \Li_{-k\alpha,0}(z)
  = \Gamma(1+k\alpha)(1-z)^{-k\alpha-1} + O(|1-z|^{-k\alpha-\epsilon})
  \]
  by Lemma~\ref{lem:litoomz}.  We shall find that the dominant terms
  in the sum in~\eqref{eq:3} are those with (i)~$k_3 = 0$, (ii)~$(k_1,
  k_2, k_3) = (k-1, 1, 0)$, and (iii)~$(k_1, k_2,  k_3) = (0, k-1,
  1)$. 

  For this paragraph, consider the case that
  \(k_1\) and \(k_2\) are both nonzero.  It follows from the
  induction hypothesis that
  \begin{align*}
    \frac{z}4 \widehat{M}_{k_1}(z)\widehat{M}_{k_2}(z) &= \frac14(1-(1-z))
    \bigl[ C_{k_1}(1-z)^{-k_1\alpha'+\tfrac12} +
    O(|1-z|^{-k_1\alpha'+\tfrac12+(2\alpha-\epsilon)})
    \bigr]\\
    &\times
    \bigl[ C_{k_2}(1-z)^{-k_2\alpha'+\tfrac12} +
    O(|1-z|^{-k_2\alpha'+\tfrac12+(2\alpha-\epsilon)})
    \bigr]\\
    &= \frac14C_{k_1}C_{k_2} (1-z)^{-(k_1+k_2)\alpha' + 1}
    + O(|1-z|^{-(k_1+k_2)\alpha'+1+(2\alpha-\epsilon)}).
  \end{align*}
  If \(k_3=0\) then the corresponding contribution to \(\widehat{R}_k(z)\) is
  \[
  \frac14\binom{k}{k_1}
  C_{k_1}C_{k_2} (1-z)^{-k\alpha' + 1}
    + O(|1-z|^{-k\alpha'+1+(2\alpha-\epsilon)}).
    \]
    If \(k_3 \ne 0\) we use Lemma~\ref{lem:omztoli} to express
    \begin{multline*}
       \frac{z}4 \widehat{M}_{k_1}(z)\widehat{M}_{k_2}(z) =
         \frac{C_{k_1}C_{k_2}}{4\Gamma((k_1+k_2)\alpha'-1)}
       \Li_{-(k_1+k_2)\alpha'+2,0}(z)\\
       +
       O(|1-z|^{-(k_1+k_2)\alpha'+1+(2\alpha-\epsilon)}) -
         \frac{C_{k_1}C_{k_2}}{4} 
       [(k_1+k_2)\alpha' <
         2]\frac{\zeta(-(k_1+k_2)\alpha'+2)}
         {\Gamma((k_1+k_2)\alpha'-1)}.
    \end{multline*}
  The corresponding contribution to \(\widehat{R}_{k}(z)\) is then
  \(\binom{k}{k_1,k_2,k_3}\) times:
  \begin{equation*}
  \frac{C_{k_1}C_{k_2}}{4\Gamma((k_1+k_2)\alpha'-1)}
  \Li_{-k\alpha'+\tfrac{k_3}2+2,0}(z)
  + \Li_{-k_3\alpha,0}(z)
  \odot O(|1-z|^{-(k_1+k_2)\alpha'+1+(2\alpha-\epsilon)}).
  \end{equation*}
  Now \(k_3 \leq k-2\) so \(-k\alpha'+\tfrac{k_3}2 + 2 < 1\).
  Hence the contribution when \(k_3 \ne 0\) is
  \[
  O(|1-z|^{-k\alpha'+\tfrac{k_3}2+1}) =
  O(|1-z|^{-k\alpha'+\tfrac32}) =
  O(|1-z|^{-k\alpha'+1+(2\alpha-\epsilon)}).
  \]

  Next we consider the case when \(k_1\) is nonzero but \(k_2=0\).  In
  this case using the induction hypothesis we see that
  \begin{align*}
    \frac{z}4 \widehat{M}_{k_1}(z)\widehat{M}_{k_2}(z) &=
    \frac{z}4 \CAT(z/4)\widehat{M}_{k_1}(z)\\
    &= \frac{1 - (1-z)^{1/2}}{2} \bigl[
    C_{k_1}(1-z)^{-k_1\alpha' + \tfrac12}\bigr]
    +
    O(|1-z|^{-k_1\alpha' + \tfrac12 +
      (2\alpha-\epsilon)})\\
    &= \frac{C_{k_1}}2 (1-z)^{-k_1\alpha'+\tfrac12} +
    O(|1-z|^{-k_1\alpha' + \tfrac12 + (2\alpha-\epsilon)}).
  \end{align*}
  Applying Lemma~\ref{lem:omztoli} to the last expression we get
  \begin{multline*}
  \frac{z}4 \widehat{M}_{k_1}(z)\widehat{M}_{k_2}(z) =
  \frac{C_{k_1}}{2\Gamma(k_1\alpha'-\tfrac12)}
  \Li_{-k_1\alpha'+\tfrac32,0}(z) \\
  + O(|1-z|^{-k_1\alpha' + \tfrac12 + (2\alpha-\epsilon)})
  - \frac{C_{k_1}}{2}[k_1\alpha'-\tfrac12 < 1]
  \frac{\zeta(-k_1\alpha'+\tfrac32)} 
    {\Gamma(k_1\alpha'-\tfrac12)}.
   \end{multline*}
   The contribution to \(\widehat{R}_{k}(z)\) is hence \(\binom{k}{k_1}\)
   times:
   \[
   \frac{C_{k_1}}{2\Gamma(k_1\alpha'-\tfrac12)}
   \Li_{-k\alpha'+\tfrac{k_3}2+\tfrac32,0}(z) +
   \Li_{-k_3\alpha,0}(z) \odot O(|1-z|^{-k_1\alpha' +
     \tfrac12 + (2\alpha-\epsilon)}).
   \]
   Using the fact that \(\alpha > 0\) and \(k_3 \leq k-1\), we conclude
   that \(-k\alpha'+\tfrac{k_3}2+\tfrac32 < 1\) so that, by
   Lemma~\ref{lem:litoomz} and
   part~(\ref{item:hadamard1}) of Theorem~\ref{thm:hadamard}, the
   contribution is
   \[
   O(|1-z|^{-k\alpha'+\tfrac{k_3}2+\tfrac12}) =
   O(|1-z|^{-k\alpha'+\tfrac32})
   \]
   where the displayed equality holds unless \(k_3=1\).  When
   \(k_3=1\) we get a corresponding contribution to $\widehat{R}_k(z)$
   of \(\binom{k}{k-1}\) times:
   \begin{equation*}
   \frac{C_{k-1}\Gamma(k\alpha'-1)}
   {2\Gamma((k-1)\alpha'-\tfrac12)}
   (1-z)^{-k\alpha'+1} +
     O(|1-z|^{-k\alpha'+1+(2\alpha-\epsilon)}),
   \end{equation*}
   since for \(k \geq 2\) we have \(k\alpha' > 1 +
   (2\alpha-\epsilon)\).  The introduction of \(\epsilon\) handles the
   case when \(k\alpha' = 1 + 2\alpha\), which would have otherwise,
   according to part~(\ref{item:hadamard3}) of Thoerem~\ref{thm:hadamard},
   introduced a logarithmic remainder.  In either case the
   remainder is $ O(|1-z|^{-k\alpha'+1+(2\alpha-\epsilon)})$.
   The case when \(k_2\) is nonzero but \(k_1=0\) is handled similarly by
   exchanging the roles of \(k_1\) and \(k_2\).
   
   The final contribution comes from the single term where both
   \(k_1\) and \(k_2\) are zero.  In this case the contribution to
   \(\widehat{R}_k(z)\) is, recalling~\eqref{eq:50},
   \begin{equation}
     \label{eq:5}
   \Li_{-k\alpha,0}(z) \odot [\frac{z}4 \CAT^2(z/4)] =
   \Li_{-k\alpha,0}(z) \odot 
   (\CAT(z/4)-1)= \Li_{-k\alpha,0}(z) \odot \CAT(z/4).
   \end{equation}
   Now, using Theorem~1 of~\cite{MR2000a:05015},
   \begin{align*}
     \CAT(z/4)             & = 2-2(1-z)^{1/2} + O(|1-z|) \\
     & = 2  +   2\frac{\zeta(3/2)}{\Gamma(-1/2)}
     - \frac{2}{\Gamma(-1/2)}\Li_{3/2,0}(z)  +
     O(|1-z|),
   \end{align*}
   so that~\eqref{eq:5} is
   \[
   -\frac2{\Gamma(-1/2)} \Li_{\tfrac32-k\alpha,0}(z) +
   O(|1-z|^{1-k\alpha}) +
   \begin{cases}
     0                     & 1-k\alpha < 0,              \\
     O(|1-z|^{-\epsilon}) & 1-k\alpha=0,                \\
     O(1)                  & 1-k\alpha > 0.
   \end{cases}
   \]
   When \(\tfrac32-k\alpha < 1\) this is
   \(O(|1-z|^{-k\alpha+\tfrac12})\); when \(\tfrac32-k\alpha \geq 1\), it
   is \(O(1)\).  In either case we get a contribution
   which is \(O(|1-z|^{-k\alpha'+1+(2\alpha-\epsilon)})\).

   Hence
   \begin{align*}
     \widehat{R}_k(z) &= \Biggl[ \sum_{\substack{k_1+k_2=k\\k_1,k_2<k}}
     \binom{k}{k_1} \frac{C_{k_1}C_{k_2}}4 + 2k \frac{C_{k-1}}2
     \frac{\Gamma(k\alpha+\tfrac{k}2-1)}{\Gamma((k-1)\alpha+\tfrac{k}2-1)}
     \Biggr] (1-z)^{-k\alpha'+1}\\
     & \qquad\qquad + O(|1-z|^{-k\alpha'+1+(2\alpha-\epsilon)})\\
     &= C_k(1-z)^{-k\alpha'+1} + O(|1-z|^{-k\alpha'
       + 1 + (2\alpha-\epsilon)}),
   \end{align*}
   with the \(C_k\)'s defined by the recurrence~\eqref{eq:10}.  Now
   using~\eqref{eq:4}, the claim follows.
 \end{proof}

\subsubsection{Large toll functions\texorpdfstring{ ($ \alpha \geq
    1/2$)}{}}
\label{sec:large-toll-functions}
When \(\alpha \geq 1/2\) there is no need to apply the centering
techinques.  Define \({\mu}_n(k) := \E{{X}_n^k}\) and
\( \bar{\mu}_n(k) := \beta_n{\mu}_n(k)/4^n \).  Let
\(\overline{M}_k(z)\) denote the 
ordinary generating function of \(\bar{\mu}_n(k)\) in $n$.  Observe that
\(\overline{M}_0(z) = \CAT(z/4)\).  As earlier, conditioning on the
key stored at the root, we get, for \(k \geq 2\),
\[
\bar{\mu}_n(k) = \frac12 \sum_{j=1}^n \frac{\beta_{n-j}}{4^{n-j}}
  \bar{\mu}_{j-1}(k) + \bar{r}_n(k),
\qquad n \geq 1,
\]
where
\begin{equation*}
\bar{r}_n(k) := 
     \frac14 \sum_{\substack{k_1+k_2+k_3=k\\
        k_1,k_2 < k}} \binom{k}{k_1,k_2,k_3} b_n^{k_3} \sum_{j=1}^n
    \bar{\mu}_{j-1}(k_1) \bar{\mu}_{n-j}(k_2),
\end{equation*}
for \(n \geq 1\) and \(\bar{r}_0(k) := \bar{\mu}_0(k) = \mu_0(k) = 0\).
Let \(\overline{R}_k(z)\) denote the ordinary generating
function of \(\bar{r}_n(k)\) in $n$.  Then
\begin{equation*}
\overline{M}_k(z) = \frac{\overline{R}_k(z)}{\sqrt{1-z}}  
\end{equation*}
and
\begin{equation}
\label{eq:12}
\overline{R}_k(z) = \sum_{\substack{k_1+k_2+k_3=k\\
    k_1,k_2<k}} \binom{k}{k_1,k_2,k_3} \bigl( B(z)^{\odot k_3}
\bigr) \odot
\bigl[ \frac{z}4 \overline{M}_{k_1}(z)\overline{M}_{k_2}(z)\bigr].
\end{equation}

We can now state the result about the asymptotics of the generating
function~\(\overline{M}_k\) when \(\alpha > 1/2\).  The case  \(\alpha=1/2\)
will be handled subsequently, in Proposition~\ref{thm:alpha=12}.
\begin{proposition}
  \label{thm:12alpha}
  Let $\epsilon > 0$ be arbitrary, and define
  \begin{equation}
    \label{eq:56}
    c :=
    \begin{cases}
      \alpha - \tfrac12      & \tfrac12 < \alpha < 1 \\
      \tfrac12 - \epsilon & \alpha = 1       \\
      \tfrac12            & \alpha > 1.
    \end{cases}
  \end{equation}
  Then the generating function $\overline{M}_k(z)$
  of~\(\bar{\mu}_n(k)\)  has the singular expansion
  \[
  \overline{M}_k(z) = C_k(1-z)^{-k(\alpha+\tfrac12)+\tfrac12} +
  O(|1-z|^{-k(\alpha+\tfrac12)+\tfrac12+c})
  \]
  for \(k \geq 1\), where the \(C_k\)'s are defined by the
  recurrence~\eqref{eq:10}.
\end{proposition}
\begin{proof}
  The proof is very similar to that of Proposition~\ref{thm:0alpha14}.
  We present a sketch.  The reader is invited to compare the cases
  enumerated below to those in the earlier proof.
  
  When \(k=1\) the claim is true by~\eqref{eq:8}.  We analyze the
  various terms in~\eqref{eq:12} for \(k \geq 2\), employing the
  notational convenience \(\alpha' := \alpha+\tfrac12\). 

  When both \(k_1\) and \(k_2\) are nonzero then the contribution to
  \(\overline{R}_k(z)\) is
  \[
  \frac14\binom{k}{k_1} C_{k_1}C_{k_2} (1-z)^{-k\alpha'+1} +
  O(|1-z|^{-k\alpha'+c+1})
  \]
  when \(k_3=0\) and is \(O(|1-z|^{-k\alpha'+c+1})\) otherwise.

  When \(k_1\) is nonzero and \(k_2=0\) the contribution to
  \(\overline{R}_k(z)\)
  is
  \[
  k \frac{C_{k-1} \Gamma(k\alpha'-1)}{2\Gamma((k-1)\alpha'-\tfrac12)}
  (1-z)^{-k\alpha'+1} + O(|1-z|^{-k\alpha'+c+1})
  \]
  when \(k_3=1\) and \(O(|1-z|^{-k\alpha'+c+1})\) otherwise.  The
  case when \(k_2\) is nonzero and \(k_1=0\) is identical.

  The final contribution comes from the single term when both \(k_1\)
  and \(k_2\) are zero.  In this case we get a contribution of
  \(O(|1-z|^{-k\alpha+\tfrac12})\) which is
  \(O(|1-z|^{-k\alpha'+c+1})\).  Adding all these contributions
  yields the desired result.
\end{proof}

The result when \(\alpha=1/2\) is as follows.  Recall that $L(z) :=
\log((1-z)^{-1})$.
\begin{proposition}
  \label{thm:alpha=12}
  Let $\alpha=1/2$.
  In the notation of Proposition~\ref{thm:12alpha},
  \[
  \overline{M}_k(z) = (1-z)^{-k+\tfrac12} \sum_{l=0}^k C_{k,l}
  L^{k-l}(z) + O(|1-z|^{-k+1-\epsilon})
  \]
  for \(k \geq 1\) and any \(\epsilon > 0\), where the $C_{k,l}$'s are
  constants.  The constant multiplying the
  lead-order term is given by
  \begin{equation}
    \label{eq:11}
    C_{k,0} = \frac{(2k-2)!}{2^{2k-2}(k-1)!\pi^{k/2}}.
  \end{equation}
\end{proposition}
\begin{proof}
  We omit the proof, referring the interested reader
  to~\cite{FK-catalan-arXiv}.
\end{proof}

\subsection{Asymptotics of moments}
\label{sec:asymptotics-moments}
For \(0 < \alpha< 1/2 \), we have seen in Proposition~\ref{thm:0alpha14}
that the generating function $\widehat{M}_k(z)$ of
\( \hat{\mu}_n(k) = \beta_n\tilde{\mu}_n(k)/4^n \) has the singular
expansion
\[
\widehat{M}_k(z) = C_k(1-z)^{-k(\alpha+\tfrac{1}2)+\tfrac12} +
O(|1-z|^{-k(\alpha+\tfrac{1}{2})+\tfrac12 + c}),
\]
where \(c := \min\{2\alpha-\epsilon,1/2\}\).
By singularity analysis~\cite{MR90m:05012},
\[
\frac{\beta_n\tilde{\mu}_n(k)}{4^n} =
C_k
\frac{n^{k(\alpha+\tfrac12)-\tfrac32}}{\Gamma(k(\alpha+\tfrac12)-\tfrac12)}
  + O(n^{k(\alpha+\tfrac{1}{2})-\tfrac32-c}).
\]
Recall that
\[
\beta_n = 
\frac{4^n}{\sqrt{\pi}n^{3/2}} \left(1 + O(\tfrac1n)\right),
\]
so that
\begin{equation}
  \label{eq:15}
  \tilde{\mu}_n(k)
  = \frac{C_k\sqrt{\pi}}{\Gamma(k(\alpha+\tfrac{1}2)-\tfrac12)}
  n^{k(\alpha+\tfrac{1}2)} +
  O(n^{k(\alpha+\tfrac{1}2) - c}).
\end{equation}

For \(\alpha > 1/2\) a similar analysis using
Proposition~\ref{thm:12alpha} yields
\begin{equation}
  \label{eq:16}
  \mu_n(k) = \frac{C_k\sqrt{\pi}}{\Gamma(k(\alpha+\tfrac{1}2)-\tfrac12)}
  n^{k(\alpha+\tfrac{1}2)} +
  O(n^{k(\alpha+\tfrac{1}2) - c}),
\end{equation}
with now $c$ as defined at~\eqref{eq:56}.  Finally, when \(\alpha=1/2\) the
asymptotics of the moments are given by
\begin{equation}
  \label{eq:17}
  \mu_n(k) =
  \left( \frac{1}{\sqrt\pi} \right)^k(n\log{n})^k +
  O(n^k(\log{n})^{k-1}).
\end{equation}

\subsection{The limiting distributions}
\label{sec:limit-distr}
In Section~\ref{sec:alpha-ne-12} we will use our moment
estimates~\eqref{eq:15} and \eqref{eq:16} with the method of moments
to derive limiting distributions for our additive functions.  The
case $\alpha=1/2$ requires a somewhat delicate analysis, which we will
present separately in Section~\ref{sec:alpha=12}.

\subsubsection{\texorpdfstring{$\alpha \ne 1/2$}{alpha not 1/2}}
\label{sec:alpha-ne-12}

We first handle the case \(0 < \alpha < 1/2\). (We assume this
restriction until just before
Proposition~\ref{thm:limit_law_alpha_ne_12}.) We have
\begin{equation}
  \label{eq:39}
  \tilde\mu_n(1) = \E\widetilde{X}_n = \E[X_n - C_0(n+1)]
  = \frac{C_1\sqrt\pi}{\Gamma(\alpha)} n^{\alpha+\tfrac12} +
  O(n^{\alpha+\tfrac12-c})
\end{equation}
with $c := \min\{2\alpha-\epsilon,1/2\}$ and
\begin{equation*}
  \tilde\mu_n(2) = \E\widetilde{X}_n^2 =
  \frac{C_2\sqrt\pi}{\Gamma(2\alpha+\tfrac12)} n^{2\alpha+1} +
  O(n^{2\alpha+1-c}).
\end{equation*}
So
\begin{equation}\label{eq:40}
  \Var{X_n} = \Var{\widetilde{X}_n} = \tilde\mu_n(2) -
  [\tilde\mu_n(1)]^2 = \sigma^2 n^{2\alpha+1} + O(n^{2\alpha+1-c}),
\end{equation}
where
\begin{equation}\label{eq:46}
  \sigma^2 := \frac{C_2\sqrt\pi}{\Gamma(2\alpha+\tfrac12)} -
  \frac{C_1^2\pi}{\Gamma^2(\alpha)}.
\end{equation}
We also have, for \(k \geq 1\),
\begin{equation}\label{eq:41}
  \E \left[ \frac{\widetilde{X}_n}{n^{\alpha+\tfrac12}} \right]^k
  = \frac{\tilde{\mu}_n(k)}{n^{k(\alpha+\tfrac12)}}
  = \frac{C_k\sqrt\pi}{\Gamma(k(\alpha+\tfrac12)-\tfrac12)} + O(n^{-c}).
\end{equation}

The following lemma provides a sufficient bound on the moments
facilitating the use of the method of moments.
\begin{lemma}
  \label{lem:Ckbound}
  Define \(\alpha' := \alpha+\tfrac12\).  There exists a constant \(A
  < \infty \)
  depending only on \(\alpha\) such that
  \begin{equation*}
    \left|\frac{C_k}{k!}\right| \leq A^k k^{\alpha'k}
  \end{equation*}
  for all \(k \geq 1\).
\end{lemma}
\begin{proof}
  The proof is fairly similar to those of
  Propositions~\ref{thm:0alpha14},~\ref{thm:12alpha} and
  Proposition~\ref{thm:shape_moments}.  We omit the details, referring
  the reader to~\cite{FK-catalan-arXiv}.
\end{proof}
It follows from Lemma~\ref{lem:Ckbound} and Stirling's approximation that
\begin{equation}
  \label{eq:19}
  \left|\frac{C_k\sqrt\pi}{k! \Gamma(k(\alpha+\tfrac12)-\tfrac12)}\right|
  \leq B^k 
\end{equation}
for large enough \(B\) depending only on $\alpha$.  Using standard
arguments~\cite[Theorem 30.1]{MR95k:60001} it follows that \(X_n\)
suitably normalized has a limiting distribution that is characterized
by its moments.  Before we state the result, we observe that the
argument presented above can be adapted with minor modifications to
treat the case \(\alpha > 1/2\), with \(\widetilde{X}_n\) replaced by
\(X_n\).  We can now state a result for $\alpha \ne 1/2$.  We will use
the notation $\stackrel{\mathcal{L}}{\to}$ to denote convergence in
law (or distribution).

\begin{proposition}\label{thm:limit_law_alpha_ne_12}
  Let $X_n$ denote the additive functional on Catalan trees
  induced by the toll sequence $(n^\alpha)_{n \geq 0}$.  Define the
  random variable $Y_n$ as follows:
  \begin{equation*}
    Y_n := 
    \begin{cases}
      \displaystyle\frac{X_n-C_0(n+1)}{n^{\alpha+\tfrac12}} & 0 <
      \alpha < 1/2,\\
      \displaystyle\frac{X_n}{n^{\alpha+\tfrac12}} & \alpha > 1/2,
    \end{cases}
  \end{equation*}
  where
  \begin{equation*}
    C_0 := \sum_{n=0}^\infty n^\alpha \frac{\beta_n}{4^n}, \qquad
    \beta_n = \frac{1}{n+1} \binom{2n}{n}.
  \end{equation*}
  Then
  \begin{equation*}
    Y_n \stackrel{\mathcal{L}}{\to} Y;
  \end{equation*}
  here $Y$ is a random variable with the unique distribution whose
  moments are
  \begin{equation}\label{eq:49}
    \E Y^k = \frac{C_k
      \sqrt\pi}{\Gamma(k(\alpha+\tfrac12)-\tfrac12)}, 
  \end{equation}
  where the $C_k$'s satisfy the recurrence
  \begin{equation*}
    C_k = \frac{1}{4} \sum_{j=1}^{k-1} \binom{k}{j} C_j C_{k-j} + k
    \frac{\Gamma(k\alpha + \tfrac{k}2 -1)}{\Gamma((k-1)\alpha +
    \tfrac{k}2 - 1)} C_{k-1}, k \geq 2; \quad C_1 =
    \frac{\Gamma(\alpha-\tfrac12)}{\sqrt\pi}.
  \end{equation*}
\end{proposition}
The case $\alpha=1/2$ is handled in Section~\ref{sec:alpha=12},
leading to Proposition~\ref{prop_alpha_12}, and a unified result for
all cases is stated as Theorem~\ref{thm:limit-dist}.

\begin{remark}\label{remark_y_alpha_properties}
We now consider some properties of the limiting random variable $Y
\equiv Y(\alpha)$ defined by its moments at~\eqref{eq:49} for $\alpha
\ne 1/2$.
\begin{enumerate}[(a)]
\item 
  When \( \alpha=1 \), setting $\Omega_k := C_k/2$ we see immediately
  that 
  \begin{equation*}
    \E Y^k =  \frac{-\Gamma(-1/2)}{\Gamma((3k-1)/2)}\Omega_k,
  \end{equation*}
  where
  \begin{equation*}
    2\Omega_k = \sum_{j=1}^{k-1} \binom{k}{j} \Omega_j \Omega_{k-j} + k(3k - 4)
    \Omega_{k-1}, \qquad \Omega_1 = \frac12.
  \end{equation*}
  Thus \(Y\) has the ubiquitous Airy distribution and we have
  recovered the limiting distribution of path length in Catalan
  trees~\cite{MR92m:60057,MR92k:60164}.  The Airy distribution arises
  in many contexts including parking allocations, hashing tables,
  trees, discrete random walks, mergesorting, etc.---see, for
  example, the introduction of~\cite{MR2002j:68115} which contains
  numerous references to the Airy distribution.
\item When $\alpha = 2$, setting $\eta := Y/\sqrt2$ and $a_{0,l} :=
  2^{2l-1} C_l$, we see that
  \begin{equation*}
    \E{\eta^l} = \frac{\sqrt\pi}{2^{(5l-2)/2} \Gamma((5l-1)/2)} a_{0,l},
  \end{equation*}
  where
  \begin{equation*}
    a_{0,l} = \frac{1}{2} \sum_{j=1}^{l-1} \binom{l}{j} a_{0,j}
    a_{0,l-j} + l(5l-4)(5l-6), \qquad a_{0,1} = 1.
  \end{equation*}
  We have thus recovered the recurrence for the moments of the
  distribution $\mathcal{L}({\eta})$, which arises in the study of the
  Wiener index  of Catalan trees~\cite[proof of Theorem~3.3 in
  Section~5]{janson:_wiener}.
\item Consider the variance $\sigma^2$ defined at~\eqref{eq:46}.
  \begin{enumerate}[(i)]
  \item   Figure~\ref{fig:variance}, plotted using
    \texttt{Mathematica}, suggests that $\sigma^2$ is positive for
    all $\alpha > 0$.  We will prove this fact in
    Theorem~\ref{thm:limit-dist}.
    \begin{figure}[htbp]
      \centering
      \includegraphics{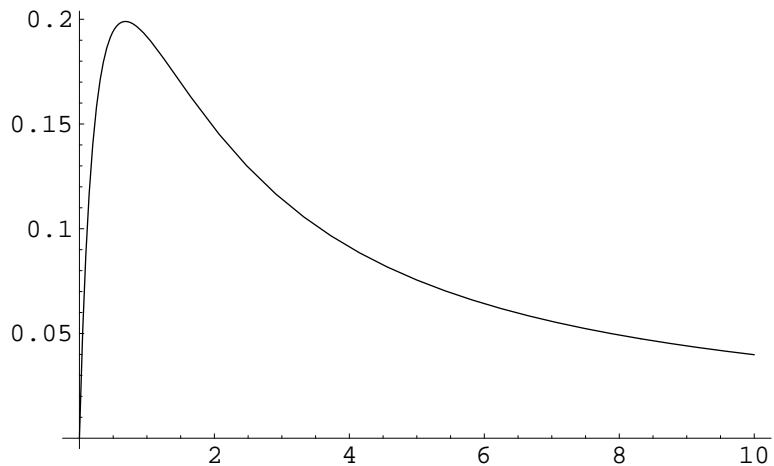}
      \caption{$\sigma^2$ of~\eqref{eq:46} as a function of $\alpha$.}
      \label{fig:variance}
    \end{figure}
    There is also numerical evidence that $\sigma^2$ is unimodal with
    $\max_{\alpha} \sigma^2(\alpha) \doteq 0.198946$ achieved at $\alpha \doteq
    0.682607$. (Here $\doteq$ denotes approximate equality.)
  \item As $\alpha \to \infty$, using Stirling's approximation one can
    show that $\sigma^2 \sim (\sqrt{2}-1)\alpha^{-1}$.
  \item As $\alpha \downarrow 0$, using a Laurent series expansion of
    $\Gamma(\alpha)$ we see that $\sigma^2 \sim 4(1-\log2) \alpha$.
  \item Though the random variable $Y(\alpha)$ has been defined only
    for $\alpha \ne 1/2$, the variance $\sigma^2$ has a limit at $\alpha=1/2$:
    \begin{equation}\label{eq:51}
      \lim_{\alpha \to 1/2} \sigma^2(\alpha) = \frac{8 \log 2}{\pi} -
      \frac{\pi}2.
    \end{equation}
  \end{enumerate}
\item Figure~\ref{fig:mc3} shows the third central moment
  $\E[Y - \E{Y}]^3$ as a function of $\alpha$. The plot
  suggests that the third central moment is positive for each $\alpha
  > 0$, which would also establish that $Y(\alpha)$ is not normal for any
  $\alpha > 0$.  However we do not know a proof of this positive
  skewness.
  [Of course, the law of
  $Y(\alpha)$ is not normal for any $\alpha > 1/2$, since its support
  is a subset of $[0,\infty)$.]
  \begin{figure}[htbp]
    \centering
    \includegraphics{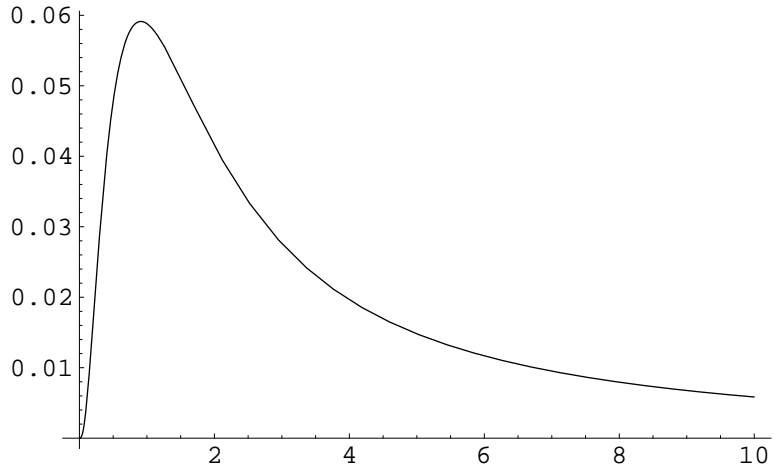}
    \caption{$\E[Y - \E{Y}]^3$ of
    Proposition~\ref{thm:limit_law_alpha_ne_12} as a function of $\alpha$.}
    \label{fig:mc3}
  \end{figure}\label{item:1}
\item When $\alpha=0$, the additive functional with toll sequence
  $(n^\alpha=1)_{n \geq 1}$ is $n$ for all trees with $n$ nodes.
  However, if one considers the random variable
  $\alpha^{-1/2}Y(\alpha)$ as $\alpha \downarrow 0$,
  using~\eqref{eq:49} and induction one can show that
  $\alpha^{-1/2}Y(\alpha)$ converges in distribution to the normal
  distribution with mean~0 and variance $4(1-\log2)$.
\item Finally, if one considers the random variable
  $\alpha^{1/2}Y(\alpha)$ as $\alpha \to \infty$, again
  using~\eqref{eq:49} and induction we find that
  $\alpha^{1/2}Y(\alpha)$ converges in distribution to the unique
  distribution with $k$th moment $\sqrt{k!}$ for $k=1,2,\ldots$.  In
  Remark~\ref{remark:sqrkfact} next, we will show that the limiting
  distribution has a bounded, infinitely smooth density on $(0,
  \infty)$.
\end{enumerate}
\end{remark}

\begin{figure}[htbp]
  \centering
  \includegraphics[width=2.4in,keepaspectratio]{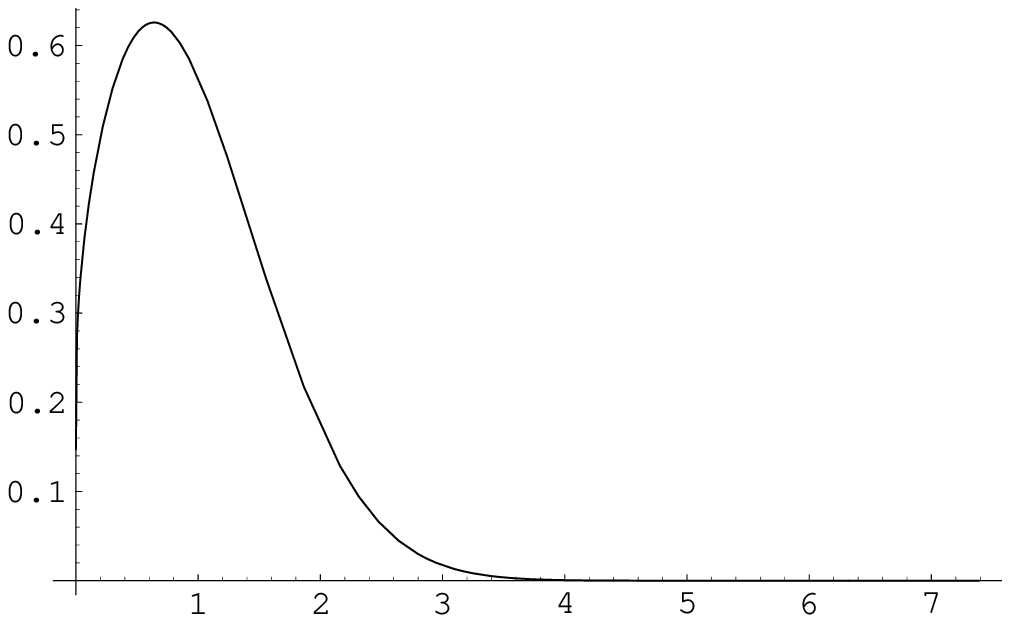}
  \includegraphics[width=2.4in,keepaspectratio]{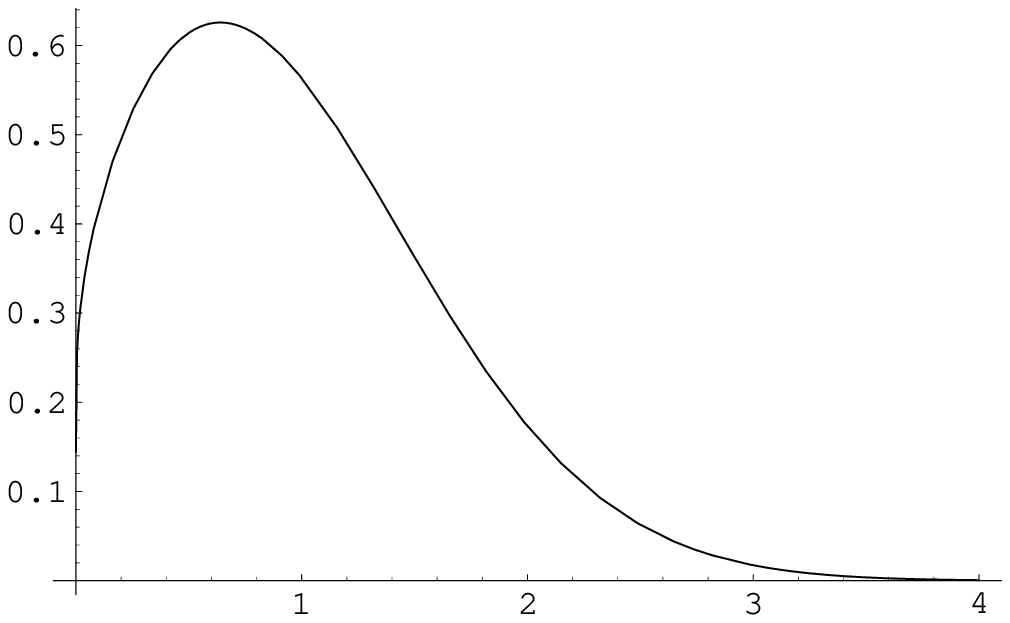}
  \caption{$f_Y$ of Remark~\ref{remark:sqrkfact}.}
  \label{fig:sqrtfactdensity}
\end{figure}

\begin{remark}
  \label{remark:sqrkfact}
  Let $Y$ be the unique distribution whose $k$th moment is $\sqrt{k!}$
  for $k=1,2,\ldots$.  Taking $Y^*$ to be an independent copy of $Y$
  and defining $X := Y Y^*$, we see immediately that $X$ is
  Exponential with unit mean.  It follows by taking logarithms that
  the distribution of $\log Y$ is a convolution square root of the
  distribution of $\log X$.  In particular, the characteristic
  function $\phi$ of $\log Y$ has square equal to $\Gamma(1 + it)$ at
  $t \in (-\infty, \infty)$; we note in passing that~$\Gamma(1 + it)$
  is the characteristic function of $-G$, where $G$ has the Gumbel
  distribution.  By exponential decay of $\Gamma(1 + it)$ as $t \to
  \pm\infty$ and standard theory (see,
  e.g.,~\cite[Chapter~XV]{MR42:5292}), $\log Y$ has an infinitely
  smooth density on $(-\infty, \infty)$, and the density and each of
  its derivatives are bounded.

  So $Y$ has an infinitely smooth density on $(0, \infty)$.  By change
  of variables, the density $f_Y$ of $Y$ satisfies
  \begin{equation*}
    f_Y(y) = \frac{f_{\log Y}(\log y)}{y}.
  \end{equation*}
  Clearly $f_Y(y)$ is bounded for $y$
  \emph{not} near~0.  (We shall drop further  consideration of
  derivatives.) To determine the behavior near~0, we need to
  know the behavior of $f_{\log Y}(\log y)/y$ as
  $y \to 0$.   Using the Fourier inversion formula, we may
  equivalently study
  \begin{equation*}
    e^x f_{\log{Y}}(-x) = \frac{1}{2\pi} \int_{-\infty}^\infty
    e^{(1+it)x} \phi(t)\, dt,
  \end{equation*}
  as $x \to \infty$.  By an application of the method of steepest
  descents [(7.2.11) in~\cite{MR89d:41049}, with $g_0=1$, $\beta =
  1/2$, $w$ the identity map, $z_0=0$, and $\alpha=0$], we get
  \begin{equation*}
    f_Y(y) \sim \frac{1}{\sqrt{\pi\log{(1/y)} }} \quad \text{as $y
      \downarrow 0$}.
  \end{equation*}
  Hence $f_Y$ is bounded everywhere.

  Using the Cauchy integral formula and simple estimates, it is easy
  to show that
  \begin{equation*}
    f_Y(y) = o( e^{-My} ) \quad \text{as $y \to \infty$}
  \end{equation*}
  for any $M < \infty$.
  Computations using the \textsc{WKB} method~\cite{MR30:3694} suggest
  \begin{equation}
    \label{eq:57}
    f_Y(y) \sim (2/\pi)^{1/4} y^{1/2} \exp(-y^2/2) \quad \text{as $y
      \to \infty$},
  \end{equation}
  in agreement with numerical calculations using
  \texttt{Mathematica}. [In fact, the right-side of~\eqref{eq:57}
  appears to be a highly accurate approximation to $f_Y(y)$ for all $y
  \geq 1$.]  Figure~\ref{fig:sqrtfactdensity} depicts the salient
  features of $f_Y$.  In particular, note the steep descent of
  $f_Y(y)$ to 0 as $y \downarrow 0$ and the quasi-Gaussian tail.
\end{remark}

\subsubsection{\texorpdfstring{$\alpha=1/2$}{alpha = 1/2}}
\label{sec:alpha=12}
For \(\alpha=1/2\), from~(\ref{eq:17}) we see immediately that
\[
\E\left[ \frac{X_n}{n\log{n}} \right]^k = \left(
  \frac1{\sqrt\pi} \right)^k + O\left( \frac{1}{\log n} \right).
\]
Thus the random variable $X_n/(n\log n)$ converges in distribution to
the degenerate random variable $1/\sqrt\pi$.  To get a nondegenerate
distribution, we carry out an analysis similar to the one that led
to~(\ref{eq:7}), getting more precise asymptotics for the mean of
$X_n$.  The refinement of~(\ref{eq:7}) that we need is the following,
whose proof we omit:
\begin{equation*}
  A(z) \odot \CAT(z/4) = \frac{1}{\sqrt\pi} (1-z)^{-1/2} L(z) +
  D_0(1-z)^{-1/2} + O(|1-z|^{\tfrac12-\epsilon}),
\end{equation*}
where
\begin{equation}
  \label{eq:58}
  D_0 = \sum_{n=1}^\infty n^{1/2} [ 4^{-n} \beta_n - \frac{1}{\sqrt{\pi}}
  n^{-3/2} ].
\end{equation}
By singularity analysis this leads to
\begin{equation}\label{eq:43}
  \E{X_n} = \frac{1}{\sqrt\pi} n \log n + D_1 n + O(n^\epsilon),
\end{equation}
where
\begin{equation}
\label{eq:59}
  D_1 = \frac{1}{\sqrt{\pi}} ( 2 \log2 + \gamma + \sqrt{\pi} D_0 ).
\end{equation}
Now analyzing the random variable $X_n -
\pi^{-1/2} n \log n$ in a manner similar to that of
Section~\ref{sec:small-toll-functions} we obtain
\begin{equation}\label{eq:44}
  \Var{[ X_n - \pi^{-1/2} n \log n ]} =
  \left(
    \frac{8}{\pi} \log 2 - \frac{\pi}{2}
  \right)
  n^2 + O(n^{\tfrac32 + \epsilon}).
\end{equation}
Using~(\ref{eq:43}) and~(\ref{eq:44}) we conclude that
\begin{equation*}
  \E
  \left[
    \frac{X_n - \pi^{-1/2} n \log n - D_1 n}{n}
  \right]
  = o(1) \end{equation*}
and
\begin{equation}\label{eq:45}
  \Var{
  \left[
    \frac{X_n - \pi^{-1/2} n \log n - D_1 n}{n}
  \right]}
  \longrightarrow \frac{8}{\pi} \log 2 - \frac{\pi}{2} = \lim_{\alpha \to 1/2}
  \sigma^2(\alpha),
\end{equation}
where $\sigma^2 \equiv \sigma^2(\alpha)$ is defined at~(\ref{eq:46})
for $\alpha \ne 1/2$. [Recall~\eqref{eq:51} of
Remark~\ref{remark_y_alpha_properties}.]

It is possible to carry out a program similar to that of
Section~\ref{sec:higher-moments} to derive asymptotics of higher order
moments using singularity analysis.  However we choose to sidestep
this arduous, albeit mechanical, computation.  Instead we will derive
the asymptotics of higher moments using a somewhat more direct
approach akin to the one employed in~\cite{MR97f:68021}.  The approach
involves approximation of sums by Riemann integrals.  To that end,
define
\begin{equation}\label{eq:47}
\widetilde{X}_n := X_n - \pi^{-1/2}(n+1) \log (n+1) -
D_1(n+1),  \quad\text{ and }\quad  \hat\mu_n(k) := \frac{\beta_n}{4^{n+1}}
\E\widetilde{X}_n^k.
\end{equation}
Note that $\widetilde{X}_0 = -D_1$, $\hat\mu_n(0) = \beta_n/4^{n+1}$,
and $\hat\mu_0(k) = (-D_1)^k/4$.  Then, in a now familiar manner, for $n
\geq 1$ we find
\begin{equation*}
  \hat\mu_n(k) = 2 \sum_{j=1}^n \frac{\beta_{j-1}}{4^j}
  \hat\mu_{n-j}(k) + \hat{r}_n(k),
\end{equation*}
where now we define
\begin{multline*}
  \hat{r}_n(k) := \sum_{\substack{k_1+k_2+k_3=k\\k_1,k_2 < k}}
  \binom{k}{k_1,k_2,k_3} \sum_{j=1}^n \hat\mu_{j-1}(k_1)
  \hat{\mu}_{n-j}(k_2)\\
  \times \left[
    \frac{1}{\sqrt\pi} ( j\log j + (n+1-j) \log (n+1-j) - (n+1) \log
  (n+1) + \sqrt\pi n^{1/2})
  \right]^{k_3}
\end{multline*}
Passing to generating functions and then back to sequences one gets,
for $n \geq 0$,
\begin{equation*}
  \hat{\mu}_n(k) = \sum_{j=0}^n (j+1) \frac{\beta_j}{4^j}
  \hat{r}_{n-j}(k). 
\end{equation*}
Using induction on $k$, we can approximate $\hat{r}_n(k)$ and
$\hat{\mu}_n(k)$ above by integrals and obtain the following result.
We omit the proof, leaving it as an exercise for the ambitious reader.
\begin{proposition}\label{prop_alpha_12}
  Let $X_n$ be the additive functional induced by the toll sequence
  $(n^{1/2})_{n \geq 1}$ on Catalan trees.  Define $\widetilde{X}_n$
  as in~{\rm (\ref{eq:47})}, with $D_1$ defined at~\eqref{eq:59} and
  $D_0$ at~\eqref{eq:58}.  Then
  \begin{equation*}
    \E
    [
      {\widetilde{X}_n}/{n}
    ]^k
    = m_k + o(1) \text{ as $n \to \infty$},
  \end{equation*}
  where $m_0=1$, $m_1=0$, and, for $k \geq 2$,
  \begin{multline}\label{eq:48}
    m_k = \frac{1}{4\sqrt\pi} \frac{\Gamma(k-1)}{\Gamma(k-\tfrac12)}\\
    \times \left[
      \sum_{\substack{k_1+k_2+k_3=k\\k_1,k_2 < k}}
    \binom{k}{k_1,k_2,k_3} m_{k_1} m_{k_2}
    \left(
      \frac{1}{\sqrt\pi}
    \right)^{k_3}
    J_{k_1,k_2,k_3}
    + 4 \sqrt\pi k m_{k-1}
    \right],
  \end{multline}
  where
  \begin{equation*}
    J_{k_1,k_2,k_3} := \int_{0}^1 x^{k_1-\tfrac32}
    (1-x)^{k_2-\tfrac32} [ x \log x + (1-x) \log (1-x) ]^{k_3} \, dx.
  \end{equation*}
  Furthermore $\widetilde{X}_n/(n+1) \stackrel{\mathcal{L}}{\to} Y$, where
  $Y$ is a random variable with the unique distribution whose
  moments are  $\E{Y^k} = m_k$, $k \geq 0$.
\end{proposition}

\subsubsection{A unified result}
\label{sec:unified-result}
The approach outlined in the preceding section can also be used for
the case $\alpha \ne 1/2$.  For completeness, we state the result for
that case here (without proof). 
\begin{proposition}\label{prop:alpha_ne_12_riemann}
  Let $X_n$ be the additive functional induced by the toll sequence
  $(n^{\alpha})_{n \geq 1}$ on Catalan trees.  Let $\alpha' := \alpha
  + \tfrac12$.  Define $\widetilde{X}_n$ as
  \begin{equation}
    \label{eq:52}
    \widetilde{X}_n :=
    \begin{cases}
      X_n - C_0(n+1) -
      \displaystyle\frac{ \Gamma(\alpha-\tfrac12
      )}{\Gamma(\alpha)}(n+1)^{\alpha'}
      & 0 < \alpha < 
      1/2, \\
      X_n - \displaystyle\frac{ \Gamma(\alpha-\tfrac12
      )}{\Gamma(\alpha)}(n+1)^{\alpha'}            & \alpha > 1/2,
    \end{cases}
  \end{equation}
  where
  \begin{equation*}
    C_0 := \sum_{n=1}^\infty n^\alpha \frac{\beta_n}{4^n}.
  \end{equation*}
  Then, for $k=0,1,2,\ldots$,
  \begin{equation*}
    \E \left[ {\widetilde{X}_n}/{n^{\alpha'}} \right]^k =
    m_k + o(1) \text{ as $n \to \infty$},
  \end{equation*}
  where $m_0=1$, $m_1=0$, and, for $k \geq 2$,
  \begin{multline}\label{eq:42}
    m_k = \frac{1}{4\sqrt\pi}
    \frac{\Gamma(k\alpha'-1)}{\Gamma(k\alpha'-\tfrac12)}\\
    \times \left[
      \sum_{\substack{k_1+k_2+k_3=k\\k_1,k_2 < k}}
      \binom{k}{k_1,k_2,k_3} m_{k_1} m_{k_2}
      \left( \frac{\Gamma(\alpha-\tfrac12)}{\Gamma(\alpha)}\right)^{k_3}
    J_{k_1,k_2,k_3}
    + 4 \sqrt\pi k m_{k-1}
    \right],
  \end{multline}
  with
  \begin{equation*}
    J_{k_1,k_2,k_3} := \int_{0}^1 x^{k_1\alpha'-\tfrac32}
    (1-x)^{k_2\alpha'-\tfrac32} [ x^{\alpha'} + (1-x)^{\alpha'} - 1
    ]^{k_3} \, dx.
  \end{equation*}
  Furthermore, $\widetilde{X}_n/n^{\alpha'}
  \stackrel{\mathcal{L}}{\to} Y$, where
  $Y$ is a random variable with the unique distribution whose moments
  are $\E{Y^k} = m_k$.
\end{proposition}

[The reader may wonder as to why we have chosen to state
Proposition~\ref{prop:alpha_ne_12_riemann} using several instances of
$n+1$, rather than $n$, in~\eqref{eq:52}.  The reason is that use of
$n+1$ is somewhat more natural in the calculations that establish the
proposition.]

In light of Propositions~\ref{thm:limit_law_alpha_ne_12},
\ref{prop_alpha_12}, and~\ref{prop:alpha_ne_12_riemann}, there are a variety of
ways to state a unified result.  We state one such version here.
\begin{theorem}
  \label{thm:limit-dist}
  Let \(X_n\) denote the additive functional induced by the toll sequence
  \((n^\alpha)_{n \geq 1}\) on Catalan trees. Then
  \begin{equation*}
    \frac{X_n - \E{X_n}}{\sqrt{\Var{X_n}}}
    \stackrel{\mathcal{L}}{\to} W,
  \end{equation*}
  where the distribution of~$W$ is described as follows:
  \bigskip\par\noindent
  (a) For $\alpha \ne 1/2$,
    \[
    W = \frac{1}{\sigma}
    \left(
      Y - \frac{C_1\sqrt\pi}{\Gamma(\alpha)}
    \right), \quad \text{ with } \quad  \sigma^2 :=
    \frac{C_2\sqrt\pi}{\Gamma(2\alpha+\tfrac12)} -
    \frac{C_1^2\pi}{\Gamma^2(\alpha)} > 0,
    \]
    where $Y$ is a random variable with the unique distribution whose
    moments are
    \begin{equation*}
      \E{Y^k} =
      \frac{C_k\sqrt\pi}{\Gamma(k(\alpha+\tfrac12)-\tfrac12)},
    \end{equation*}
    and the \(C_k\)'s satisfy the recurrence~{\rm (\ref{eq:10})}.
    \bigskip\par\noindent
    (b) For $\alpha=1/2$,
    \begin{equation*}
      W = \frac{Y}{\sigma}, \quad \text{ with } \quad \sigma^2 :=
      \frac{8}{\pi} \log 2 - \frac{\pi}{2},
    \end{equation*}
    where $Y$ is a random variable with the unique distribution whose
    moments $m_k = \E{Y^k}$ are given by~{\rm (\ref{eq:48})}.
\end{theorem}
\begin{proof}
  Define
  \begin{equation*}
    W_n := \frac{X_n - \E{X_n}}{\sqrt{\Var{X_n}}}
  \end{equation*}
  \bigskip\par\noindent
  (a) Consider first the case $\alpha < 1/2$ and let $\alpha' := \alpha +
    \tfrac12$.  By~\eqref{eq:39},
    \begin{equation}\label{eq:53}
      \E{X_n} = C_0(n+1) + \frac{C_1\sqrt\pi}{\Gamma(\alpha)}
      n^{\alpha'} + o(n^{\alpha'}).
    \end{equation}
    Since $\widetilde{X}_n$ defined at~\eqref{eq:52} and $X_n$ differ
    by a deterministic amount, $\Var{X_n}=\Var{\widetilde{X}_n}$. Now by
    Proposition~\ref{prop:alpha_ne_12_riemann},
    \begin{equation}\label{eq:54}
      \Var{\widetilde{X}_n} = \E{\widetilde{X}_n^2} -
      (\E{\widetilde{X}_n})^2 = (m_2 + o(1))n^{2\alpha'} - (m_1^2 +
      o(1))n^{2\alpha'} = (m_2 + o(1))n^{2\alpha'}.
    \end{equation}
    So $\sigma^2$ equals $m_2$ 
    defined at~(\ref{eq:42}), namely,
    \begin{equation*}
      \frac{1}{4\sqrt\pi}\frac{\Gamma(2\alpha'-1)}{\Gamma(2\alpha'-\tfrac12)}
      \left(\frac{\Gamma(\alpha-\tfrac12)}{\Gamma(\alpha)}\right)^2 J_{0,0,2}.
    \end{equation*}
    Thus to show $\sigma^2 > 0$ it is enough to show that $J_{0,0,2} > 0$.  But
    \begin{equation*}
      J_{0,0,2} = \int_{0}^1 x^{-3/2} (1-x)^{-3/2} [ x^{\alpha'} +
      (1-x)^{\alpha'} - 1]^2 \,dx,
    \end{equation*}
    which is clearly positive.  Using~\eqref{eq:53} and~\eqref{eq:54},
    \begin{equation*}
      W_n = \frac{X_n - C_0(n+1) - \frac{C_1\sqrt\pi}{\Gamma(\alpha)}
      n^{\alpha'} + o(n^{\alpha'})}{ (1 + o(1))\sigma
      n^{\alpha'}},
    \end{equation*}
    so, by Proposition~\ref{thm:limit_law_alpha_ne_12} and Slutsky's
    theorem~\cite[Theorem~25.4]{MR95k:60001}, the claim
    follows.
    
    The case $\alpha > 1/2$ follows similarly.
    \bigskip\par\noindent
  (b) When $\alpha=1/2$, 
  \begin{equation*}
    \E{X_n} = \frac{1}{\sqrt\pi} n\log{n} + D_1 n + o(n)
  \end{equation*}
  by~\eqref{eq:43} and
  \begin{equation*}
    \Var{X_n} =
    \left(
      \frac{8}{\pi}\log2 - \frac{\pi}{2} + o(1)
    \right)n^2
  \end{equation*}
  by~\eqref{eq:45}. The claim then follows easily from
    Proposition~\ref{prop_alpha_12} and Slutsky's
    theorem.
\end{proof}

\section{The shape functional}
\label{sec:shape-functional}
We now turn our attention to the shape functional for Catalan trees.
The shape functional is the cost induced by the toll function
$b_n \equiv \log{n}$, $n \geq 1$.  For background and results on the shape
functional, we refer the reader to~\cite{MR97f:68021} and~\cite{MR99j:05171}.

In the sequel we will improve on the mean and variance estimates
obtained in~\cite{MR97f:68021} and derive a central limit theorem for
the shape functional for Catalan trees.  The technique employed is
singularity analysis followed by the method of moments.

\subsection{Mean}
\label{sec:shape_function_mean}

We use the notation and techniques of
Section~\ref{sec:asympotics-mean} again.  Observe that now $B(z) =
\Li_{0,1}(z)$ and~\eqref{eq:21} gives the singular expansion
\begin{multline*}
  \CAT(z/4) = 2 - \frac2{\Gamma(-1/2)}[ \Li_{3/2,0}(z) -
  \zeta(3/2)]\\
  + 2\left(1-\frac{\zeta(1/2)}{\Gamma(-1/2)} \right)(1-z) +
  O(|1-z|^{3/2}).
\end{multline*}
So
\begin{equation*}
  B(z)\odot\CAT(z/4) = -\frac2{\Gamma(-1/2)} \Li_{3/2,1}(z)
  + \bar{c} + \bar{\bar{c}}(1-z) +  O(|1-z|^{\tfrac32-\epsilon}),
\end{equation*}
where \( \bar{c} \) and \( \bar{\bar{c}} \) denote unspecified (possibly 0)
constants.  The constant term in the singular expansion of
\( B(z)\odot\CAT(z/4) \) is already known to be
\begin{equation*}
  C_0 = B(z)\odot\CAT(z/4) \Bigr\rvert_{z=1} = \sum_{n=1}^\infty
  (\log{n}) \frac{\beta_n}{4^n}.
\end{equation*}
Now using the singular expansion of \( \Li_{3/2,1}(z) \), we get
\begin{equation*}
  B(z) \odot \CAT(z/4) =  C_0 - 2(1-z)^{1/2}L(z) -
  2(2(1-\log(2))-\gamma)(1-z)^{1/2} + O(|1-z|),
\end{equation*}
so that
\begin{equation}
  \label{eq:22}
  A(z)\odot\CAT(z/4) = C_0(1-z)^{-1/2} - 2L(z) -
  2(2(1-\log2)-\gamma) + O(|1-z|^{1/2}).
\end{equation}
Using singularity analysis and the asymptotics of the Catalan numbers
we get that the mean \( a_n \) of the shape functional  is given by
\begin{equation}
  \label{eq:24}
  a_n = C_0(n+1) - 2\sqrt\pi n^{1/2} + O(1),
\end{equation}
which agrees with the estimate in Theorem~3.1 of~\cite{MR97f:68021}
and improves the remainder estimate.

\subsection{Second moment and variance}
\label{sec:shape_functional_variance}

We now derive the asymptotics of the approximately centered second
moment and the variance of the shape functional.  These estimates will
serve as the basis for the induction to follow.  We will use the
notation of Section~\ref{sec:small-toll-functions}, centering the cost
function as before by $C_0(n+1)$. 

It is clear from~\eqref{eq:22} that
\begin{equation}\label{eq:38}
  \widehat{M}_1(z) = -2L(z) - 2(2(1-\log2)-\gamma) +
  O(|1-z|^{1/2}),
\end{equation}
and~\eqref{eq:3} with \( k=2 \) gives us, recalling~\eqref{eq:50},
\begin{equation}
  \label{eq:34}
  \widehat{R}_2(z) = C_0^2 + \CAT(z/4)\odot\Li_{0,2}(z) +  4\Li_{0,1}(z)
  \odot [\frac{z}4 \CAT(z/4)
  \widehat{M}_1(z)] + \frac{z}2\widehat{M}_1^2(z).
\end{equation}
We analyze each of the terms in this sum.  For the last term, observe
that $z/2 \to 1/2$ as $z \to 1$, so that
\begin{equation*}
  \frac{z}2 \widehat{M}_1^2(z) = 2L^2(z) + 4(2(1-\log2)-\gamma)L(z) +
  2(2(1-\log2)-\gamma)^2 + O(|1-z|^{\tfrac12-\epsilon}),
\end{equation*}
the \( \epsilon \) introduced to avoid logarithmic remainders.  The first
term is easily seen to be
\begin{equation*}
 \CAT(z/4) \odot \Li_{0,2}(z) = K + O(|1-z|^{\tfrac12-\epsilon}),
\end{equation*}
where
\begin{equation*}
  K := \sum_{n=1}^\infty (\log{n})^2 \frac{\beta_n}{4^n}.
\end{equation*}
For the middle term, first observe that
\begin{equation*}
  \frac{z}4 \CAT(z/4)\widehat{M}_1(z) = -L(z) - (2(1-\log2)-\gamma) +
  (1-z)^{1/2}L(z) + O(|1-z|^{1/2})
\end{equation*}
and that \( L(z) = \Li_{1,0}(z) \). Thus the third term on the right
in~\eqref{eq:34} is 4 times:
\begin{equation*}
  -\Li_{1,1}(z) + \bar{c} + O(|1-z|^{\tfrac12-2\epsilon}) = -\frac12
   L^2(z) + \gamma L(z) + \bar{c} + O(|1-z|^{\tfrac12-\epsilon}).
\end{equation*}
[The singular expansion for \( \Li_{1,1}(z) \)
was obtained using the results at the bottom of p.~379
in~\cite{MR2000a:05015}.  We state it here for the reader's
convenience:  
\begin{equation*}
  \Li_{1,1}(z) = \frac12 L^2(z) - \gamma L(z) + \bar{c} + O(|1-z|),
\end{equation*}
where $\bar{c}$ is again an unspecified constant.]  Hence
\begin{equation*}
  \widehat{R}_2(z) = 8(1-\log2)L(z) + \bar{c} +
  O(|1-z|^{\tfrac12-\epsilon}), 
\end{equation*}
which leads to
\begin{equation}
  \label{eq:23}
  \widehat{M}_2(z) = 8(1-\log2)(1-z)^{-1/2}L(z) +
  \bar{c}(1-z)^{-1/2} + O(|1-z|^{-\epsilon}).
\end{equation}
We draw the attention of the reader to the cancellation
of the ostensible lead-order term \( L^2(z) \).  This kind of
cancellation will appear again in the next section when we deal with
higher moments.

Now using singularity analysis and estimates for the Catalan numbers
we get
\begin{equation}
  \label{eq:25}
  \tilde\mu_n(2) = 8(1-\log2)n\log{n} + \bar{c}n + O(n^{\tfrac12+\epsilon}).
\end{equation}
Using~\eqref{eq:24},
\begin{equation*}
  \Var{X_n} = \tilde\mu_n(2) - \tilde\mu_n(1)^2 =
  8(1-\log2)n\log{n} + \bar{c}n + O(n^{\tfrac12+\epsilon}),
\end{equation*}
which agrees with Theorem~3.1 of~\cite{MR97f:68021} (after a
correction pointed out in~\cite{MR99j:05171}) and improves the
remainder estimate.  In our subsequent analysis we will not need to
evaluate the unspecified constant $\bar{c}$.

\subsection{Higher moments}
\label{sec:shape_function_higher-moments}

We now turn our attention to deriving the asymptotics of higher
moments of the shape functional.  The main result is as follows.
\begin{proposition}
  \label{thm:shape_moments}
  Define \(\widetilde{X}_n := X_n - C_0(n+1)\), with \(X_0 := 0\);
  \(\tilde{\mu}_n(k) := \E{\widetilde{X}_n^k} \), with
  $\tilde{\mu}_n(0) = 1$ for
  all $n \geq 0$; and \(\hat{\mu}_n(k)
  := \beta_n\tilde{\mu}_n(k)/4^n \).  Let \(\widehat{M}_k(z)\) denote
  the ordinary 
  generating function of \(\hat{\mu}_n(k)\) in the argument $n$.  For
  \( k \geq 2 \), \( \widehat{M}_k(z) \) has the singular expansion
  \begin{equation*}
    \widehat{M}_k(z) = (1-z)^{-\tfrac{k-1}2}
    \sum_{j=0}^{\lfloor {k}/2 \rfloor} C_{k,j}
    L^{\lfloor {k}/2 \rfloor-j}(z) +
    O(|1-z|^{-\tfrac{k}2+1-\epsilon}), 
  \end{equation*}
  with
  \begin{equation*}
    C_{2l,0} = \frac14 \sum_{j=1}^{l-1} \binom{2l}{2j}
    C_{2j,0}C_{2l-2j,0}, \qquad C_{2,0} = 8(1-\log2).
  \end{equation*}
\end{proposition}
\begin{proof}
  The proof is by induction.  For \( k=2 \) the claim is true
  by~\eqref{eq:23}.  We note that the claim is \emph{not} true for
  \( k=1 \). Instead, recalling~(\ref{eq:38}),
  \begin{equation}
    \label{eq:28}
    \widehat{M}_1(z) = -2L(z) -
    2(2(1-\log2)-\gamma) + O(|1-z|^{1/2}).
  \end{equation}
  For the induction step, let \( k \geq 3 \).  We will first get the
  asymptotics of \( \widehat{R}_k(z) \) defined at~\eqref{eq:3} with
  $B(z) = \Li_{0,1}(z)$.  In order to do that we will obtain the
  asymptotics of each term in the defining sum.  We remind the reader
  that we are only interested in the form of the asymptotic expansion
  of \( \widehat{R}_k(z) \) and the coefficient of the lead-order term
  when \( k \) is even.  This allows us to ``define away'' all other
  constants, their determination delayed to the time when the need
  arises.

  For this paragraph suppose that \( k_1 \geq2 \) and \( k_2 \geq 2 \).  Then
  by the induction hypothesis
  \begin{multline}
    \label{eq:26}
    \frac{z}4 \widehat{M}_{k_1}(z)\widehat{M}_{k_2}(z) = \frac14
    (1-z)^{-\tfrac{k_1+k_2}2+1}
    \sum_{l=0}^{\lfloor {k_1}/2 \rfloor +
    \lfloor {k_2}/2 \rfloor} A_{k_1,k_2,l}
    L^{\lfloor {k_1}/2 \rfloor +  \lfloor {k_2}/2
      \rfloor - l}(z)\\
    {}+ O(|1-z|^{-\tfrac{k_1+k_2}2+\tfrac32-\epsilon}),
  \end{multline}
  where \( A_{k_1,k_2,0}= C_{k_1,0}C_{k_2,0} \).  (a) If \( k_3=0 \)
  then \( k_1+k_2=k \) and the corresponding contribution to \(
  \widehat{R}_k(z) \) is given by
  \begin{multline}
    \label{eq:27}
    \frac14 \binom{k}{k_1} (1-z)^{-\tfrac{k}2+1}\\
    \times
    \sum_{l=0}^{\lfloor {k_1}/2 \rfloor +
      \lfloor ({k-k_1})/2 \rfloor} A_{k_1,k-k_1,l}
    L^{\lfloor {k_1}/2 \rfloor +  \lfloor ({k-k_1})/2 \rfloor -
      l}(z) + O(|1-z|^{-\tfrac{k}2+\tfrac32-\epsilon}).
  \end{multline}
  Observe that if \( k \) is even and \( k_1 \) is odd the highest power of
  \( L(z) \) in~\eqref{eq:27} is \( \lfloor {k}/2 \rfloor-1 \).  In all
  other cases the the highest power of \( L(z) \) in~\eqref{eq:27} is
  \( \lfloor {k}/2 \rfloor \).  (b) If \( k_3 \ne 0 \) then we  use
  Lemma~\ref{lem:omztoli} to express~\eqref{eq:26} as a linear
  combination of
  \begin{equation*}
    \left\{
      \Li_{-\tfrac{k_1+k_2}2+2,l}(z) \right\}_{l=0}^{\lfloor {k_1}/2 \rfloor
      + \lfloor {k_2}/2 \rfloor}
  \end{equation*}
  with a remainder that is
  $O(|1-z|^{-\tfrac{k_1+k_2}2+\tfrac32-\epsilon})$.  When we take the 
  Hadamard product of such a term with \( \Li_{0,k_3}(z) \) we will get a
  linear combination of
    \begin{equation*}
    \left\{
      \Li_{-\tfrac{k_1+k_2}2+2,l+k_3}(z) \right\}_{l=0}^{\lfloor
      {k_1}/2 \rfloor + \lfloor {k_2}/2 \rfloor}
  \end{equation*}
  and a smaller remainder.  Such terms are all
  \( O(|1-z|^{-\tfrac{k_1+k_2}2+1-\epsilon}) \), so that the
  contribution is
  \( O(|1-z|^{-\tfrac{k}2+\tfrac32-\epsilon}) \).

  Next, consider the case when \( k_1=1 \) and \( k_2 \geq 2 \).  Using the
  induction hypothesis and~\eqref{eq:28} we get
  \begin{equation}
    \label{eq:29}
    \begin{split}
    \frac{z}4 \widehat{M}_{k_1}(z) \widehat{M}_{k_2}(z) = -\frac12
    (1-z)^{-\tfrac{k_2-1}2} \sum_{j=0}^{\lfloor {k_2}/2 \rfloor+1}
    B_{k_2,j} L^{\left\lfloor \tfrac{k_2}2 \right\rfloor + 1 - j}(z) \\
       {}+ O(|1-z|^{-\tfrac{k_2}2+1-2\epsilon}),
    \end{split}
  \end{equation}
  with \( B_{k_2,0} = C_{k_2,0} \).  (a) If \( k_3=0 \) then \(
  k_2=k-1 \) and the 
  corresponding contribution to \( \widehat{R}_k(z) \) is given by
  \begin{equation}
    \label{eq:30}
    -\frac{k}{2} (1-z)^{-\tfrac{k}2+1}
     \sum_{j=0}^{\lfloor ({k-1})/{2} \rfloor+1} B_{k-1,j}
     L^{\left\lfloor \tfrac{k-1}2 \right\rfloor+1-j}(z) +
     O(|1-z|^{-\tfrac{k}2+\tfrac32-2\epsilon}).
   \end{equation}
   (b) If \( k_3 \ne 0 \) then Lemma~\ref{lem:omztoli} can be used once again
   to express~\eqref{eq:29} in terms of generalized polylogarithms,
   whence an argument similar to that at the end of the preceding paragraph
   yields that the contributions to $\widehat{R}(z)$ from such terms is
   \( O(|1-z|^{-\tfrac{k_2-1}2-\epsilon}) \), which is
   \( O(|1-z|^{-\tfrac{k}2+\tfrac32-\epsilon}) \).   The case when \(k_1
   \geq 2\)  and  \(k_2 =1\) is handled symmetrically.

   When \( k_1=k_2=1 \) then \(
   (z/4)\widehat{M}_{k_1}(z)\widehat{M}_{k_2}(z) \) is 
   \( O(|1-z|^{-\epsilon}) \) and when one takes the Hadamard product of
   this term with \( \Li_{0,k_3}(z) \) the contribution will be
   \( O(|1-z|^{-2\epsilon}) \).

   Now consider the case when \( k_1=0 \) and \( k_2 \geq 2 \).  Since
   \( \widehat{M}_0(z) = \CAT(z/4) \), we have
   \begin{equation}
     \label{eq:31}
     \frac{z}4 \widehat{M}_{k_1}(z)\widehat{M}_{k_2}(z) = \frac12
     (1-z)^{-\tfrac{k_2-1}2} \sum_{j=0}^{\lfloor {k_2}/2 \rfloor}
     C_{k_2,j} L^{\lfloor {k_2}/2 \rfloor-j}(z) +
     O(|1-z|^{-\tfrac{k_2}2+1-\epsilon}).
   \end{equation}
   By Lemma~\ref{lem:omztoli} this can be expressed as a linear
   combination of
   \begin{equation*}
     \left\{ \Li_{-\tfrac{k_2-1}2+1,j}(z)
     \right\}_{j=0}^{\lfloor {k_2}/2 \rfloor}
   \end{equation*}
   with a \( O(|1-z|^{-\tfrac{k_2}2+1-\epsilon}) \) remainder.  When we
   take the Hadamard product of such a term with \( \Li_{0,k_3}(z) \) we
   will get a linear combination, call it $S(z)$, of
   \begin{equation*}
     \left\{ \Li_{-\tfrac{k_2-1}2+1,j+k_3}(z)
     \right\}_{j=0}^{\lfloor {k_2}/2 \rfloor}
   \end{equation*}
   with a remainder of \( O(|1-z|^{-\tfrac{k_2}2 + 1 - 2\epsilon})
   \), which is  \( O(|1-z|^{-\tfrac{k}2 + \tfrac32 - 2\epsilon}) \)
   unless \( k_2=k-1 \).  When 
   \( k_2 =  k-1 \), by Lemma~\ref{lem:omztoli} the constant multiplying
   the lead-order term 
   \( \Li_{-\tfrac{k}{2}+2,\lfloor\tfrac{k-1}{2} \rfloor + 1}(z) \) in $S(z)$ is
   $\frac{C_{k-1,0}}2
   \mu_0^{(-\tfrac{k}{2}+2,\lfloor\tfrac{k-1}{2}\rfloor)}$.  When we
   take the Hadamard product of this term with \( \Li_{0,k_3}(z) \) we get
   a lead-order term of
   \begin{equation*}
   \frac{C_{k-1,0}}2
   \mu_0^{(-\tfrac{k}{2}+2,\lfloor\tfrac{k-1}{2}\rfloor)}
   \Li_{-\tfrac{k}{2}+2,\lfloor\tfrac{k-1}{2}\rfloor+1}(z).
   \end{equation*}
   Now we use Lemma~\ref{lem:litoomz} and the observation that
   \( \lambda_0^{(\alpha,r)}\mu_0^{(\alpha,s)}=1 \) to conclude that the
   contribution to $\widehat{R}_k(z)$ from the term with \( k_1=0 \)
   and \( k_2=k-1 \) is 
   \begin{equation}
     \label{eq:32}
     \frac{k}2 (1-z)^{-\tfrac{k}2+1}
     \sum_{j=0}^{\lfloor\frac{k-1}{2}\rfloor+1} D_{k,j}
     L^{\lfloor\frac{k-1}{2}\rfloor+1-j}(z)
     + O(|1-z|^{-\tfrac{k}2+\tfrac{3}2-\epsilon}),
   \end{equation}
   with \( D_{k,0}=C_{k-1,0} \).  Notice that the lead order from this
   contribution is precisely that from~\eqref{eq:30} but with opposite
   sign; thus the two contributions cancel each other to lead order.
   The case $k_2 = 0$ and $k_1 \geq 2$ is handled symmetrically.
   
   The last two cases are \( k_1=0 \), \( k_2=1 \) (or vice-versa)
   and \( k_1=k_2=0 \).  The contribution from these cases can be
   easily seen to be \( O(|1-z|^{-\tfrac{k}2+\tfrac32-2\epsilon}) \).
   
   We can now deduce the asymptotic behavior of \( \widehat{R}_k(z)
   \).  The three contributions are~\eqref{eq:27}, \eqref{eq:30},
   and~\eqref{eq:32}, with only~\eqref{eq:27} (in net) contributing a
   term of the form \( (1-z)^{-\tfrac{k}{2}+1} L^{\lfloor {k}/{2}
     \rfloor}(z) \) when \(k\) is even.  The coefficient of this term
   when \( k \) is even is given by
   \begin{equation*}
     \frac14 \sum_{\substack{0 < k_1 < k\\k_1\text{ even}}}
     \binom{k}{k_1} C_{k_1,0}C_{k_2,0}.
   \end{equation*}
   Finally we can sum up the rest of the contribution, define
   \( C_{k,j} \) appropriately and use~\eqref{eq:4} to claim the result.
\end{proof}

\subsection{A central limit theorem}
\label{sec:shape_functional_centr-limit-theor}

Proposition~\ref{thm:shape_moments} and singularity analysis allows us to
get the asymptotics of the moments of the ``approximately centered'' shape
functional.  Using arguments identical to those in
Section~\ref{sec:asymptotics-moments} it is clear that for \( k \geq 2 \)
\begin{equation*}
  \label{eq:33}
  \tilde{\mu}_n(k) = \frac{C_{k,0}\sqrt\pi}{\Gamma(\tfrac{k-1}2)}
  n^{{k}/{2}} [\log{n}]^{\lfloor{k}/{2}\rfloor} +
  O(n^{{k}/{2}} [\log{n}]^{\lfloor{k}/{2}\rfloor-1}).
\end{equation*}
This and the asymptotics of the mean derived in
Section~\ref{sec:shape_function_mean} give us, for \( k \geq 1 \),
\begin{equation*}
  E\left[ \frac{\tilde{X}_n}{\sqrt{n\log{n}}} \right]^{2k} \to
  \frac{C_{2k,0}\sqrt\pi}{\Gamma(k-\tfrac12)}, \qquad E\left[
  \frac{\tilde{X}_n}{\sqrt{n\log{n}}} \right]^{2k-1}  = o(1)
\end{equation*}
as $n \to \infty$.
The recurrence for \( C_{2k,0} \) can be solved easily to yield, for $k
\geq 1$,
\begin{equation*}
  C_{2k,0} = \frac{(2k)!(2k-2)!}{2^k2^{2k-2}k!(k-1)!} \sigma^{2k},
\end{equation*}
where \( \sigma^2 := 8(1-\log2) \). Then using the identity
\begin{equation*}
  \frac{\Gamma(k-\tfrac12)}{\sqrt\pi} = 
  \left[2^{2k-2} \frac{(k-1)!}{(2k-2)!}\right]^{-1}
\end{equation*}
we get
\begin{equation*}
  \frac{C_{2k,0}\sqrt\pi}{\Gamma(k-\tfrac12)} = \frac{(2k)!}{2^k
  k!}\sigma^{2k}.
\end{equation*}
It is clear now that both the ``approximately centered'' and the
normalized shape functional are asymptotically normal.
\begin{theorem}
  \label{thm:shape_clt}
  Let \( X_n \) denote the shape functional, induced by the toll
  sequence $(\log{n})_{n \geq 1}$, for Catalan trees.  Then
  \begin{equation*}
    \frac{X_n-C_0(n+1)}{\sqrt{n\log{n}}} \stackrel{\mathcal{L}}{\to}
    \mathcal{N}(0,\sigma^2) \quad\text{ and }\quad
    \frac{X_n - \E{X_n} }{\sqrt{ \Var{X_n} }}
    \stackrel{\mathcal{L}}{\to}  \mathcal{N}(0,1),
  \end{equation*}
  where
  \begin{equation*}
  C_0 := \sum_{n=1}^\infty (\log{n}) \frac{\beta_n}{4^n}, \qquad
  \beta_n = \frac1{n+1}\binom{2n}{n},
  \end{equation*}
  and \( \sigma^2 := 8(1-\log2) \).
\end{theorem}
Concerning numerical evaluation of the constant~$C_0$, see the end of
Section~5.2 in~\cite{FFK}.

\section{Sufficient conditions for asymptotic normality}
\label{sec:suff-cond-asympt}

In this speculative final section we briefly examine the behavior of a
general additive functional $X_n$ induced by a given ``small'' toll
sequence $(b_n)$.  We have seen evidence
[Remark~\ref{remark_y_alpha_properties}(\ref{item:1})] that if
$(b_n)$ is the ``large'' toll sequence $n^{\alpha}$ for any fixed $\alpha
>  0$, then the limiting behavior is non-normal.  When $b_n = \log n$ (or
$b_n = n^\alpha$ and $\alpha \downarrow 0$), the (limiting) random
variable is normal.  Where is the interface between normal and non-normal
asymptotics?  We have carried out arguments similar to those leading to
Propositions~\ref{prop_alpha_12} and~\ref{prop:alpha_ne_12_riemann} (see
also~\cite{MR97f:68021}) that suggest a sufficient condition for asymptotic
normality, but our ``proof'' is somewhat heuristic, and further technical
conditions on $(b_n)$ may be required.  Nevertheless, to inspire further
work, we present our preliminary indications.

We assume that $b_n \equiv b(n)$, where $b(\cdot)$ is a function of a
nonnegative real argument.  Suppose that $x^{-3/2} b(x)$ is (ultimately)
nonincreasing and that $x b'(x)$ is slowly varying at infinity.  Then
\begin{equation*}
   \E{X_n} = C_0 (n+1) - (1 + o(1)) 2 \sqrt{\pi} n^{3/2} b'(n),
\end{equation*}
where
\begin{equation*}
   C_0 = \sum_{n=1}^\infty b_n \frac{\beta_n}{4^n}.
\end{equation*}
Furthermore,
\begin{equation*}
   \Var{X_n} \sim 8 (1 - \log 2) [n b'(n)]^2 n\log n,
\end{equation*}
and
\begin{equation*}
   \frac{X_n - C_0 (n + 1)}{n b'(n) \sqrt{n \log n}}
   \stackrel{\mathcal{L}}{\to} \mathcal{N}(0,\sigma^2), \text{ where }
   \sigma^2 = 8(1 - \log 2).
\end{equation*}
This asymptotic normality can also be stated in the form
\begin{equation*}
   \frac{X_n - \E{X_n} }{\sqrt{\Var{X_n}}}
   \stackrel{\mathcal{L}}{\to} \mathcal{N}(0,1).
\end{equation*}

\medskip
\noindent
\textbf{Acknowledgments.}  We thank two anonymous referees for helpful
comments.

\bibliographystyle{habbrv}
\bibliography{msn,leftovers}

\end{document}